\documentclass[11pt]{elsarticle}
\usepackage{color}
\usepackage{graphicx,multirow}
\usepackage{amssymb,mathtools}
\usepackage{epstopdf}
\usepackage{bm}
\usepackage{float}
\usepackage[table]{xcolor}
\usepackage{subcaption}
\usepackage{mathtools}
\usepackage{float}																					
\usepackage[margin=1in,papersize={8.5in,11in}]{geometry}
\usepackage{times}
\usepackage{tabularx}
\usepackage[normalem]{ulem}

\newtheorem{remark}{\bf Remark}[section]

\DeclareMathOperator*{\argmin}{\arg\!\min}

\def\vs{\vspace{0.2cm}}

\journal{arXiv}

\begin{document}
\begin{frontmatter}
\title{A new scalable algorithm for computational optimal control under uncertainty}

\author[ucsc]{Panos Lambrianides}
\author[ucsc]{Qi Gong}
\author[ucsc]{Daniele Venturi\corref{correspondingAuthor}}
\address[ucsc]{Department of Applied Mathematics\\University of California Santa Cruz\\ Santa Cruz, CA 95064}
\cortext[correspondingAuthor]{Corresponding author}
\ead{venturi@ucsc.edu}

\begin{abstract} 
We address the design and synthesis of optimal 
control strategies for high-dimensional stochastic 
dynamical systems.
Such systems may be deterministic nonlinear systems 
evolving from random initial states, or systems driven
by random parameters or processes.
The objective is to provide a validated new computational
capability for optimal control which will be achieved more 
efficiently than current state-of-the-art methods. 
The new framework utilizes direct single or multi-shooting 
discretization, and is based on efficient vectorized gradient 
computation with adaptable memory management. 
The algorithm is demonstrated to be scalable to high 
dimensional  nonlinear control systems with 
random initial condition and unknown parameters.	
\end{abstract}

\end{frontmatter}

\section{Introduction}
\begin{figure}[!t]
	\centerline{ 
	\includegraphics[height=4.5cm]{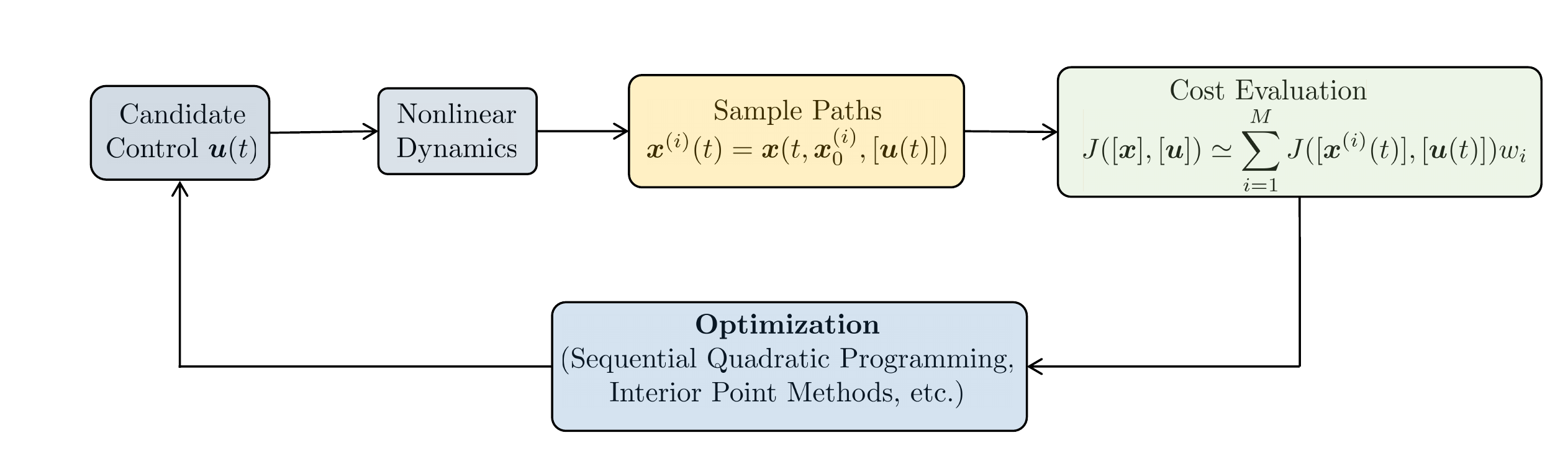}}
	\caption{Sketch of the proposed data-driven optimal
control architecture. The new control algorithm is 
based on interior point optimization methods, common 
sub-expression elimination, and exact gradient information 
obtained with automatic differentiation and computational graphs.}
	\label{fig:sketch}
\end{figure}

Designing optimal control strategies for high-dimensional 
stochastic dynamical systems is critical in many engineering 
applications, such as search of unknown targets with 
autonomous vehicles, path planning of 
heterogeneous agents, and quantum control
\cite{Foraker.16I,Automatica_search,waltonoptimal,Li1879,wang2018free}.
Such systems may be modeled as nonlinear control 
systems of the form 
\begin{equation}
\dot {\bm x}  =  \bm f(\bm x,\bm u), \qquad 
\bm x(0)=\bm x_0(\omega),\qquad t\in [0,t_f], 
\label{eq:dynamics}
\end{equation}
where $\bm x(t)\in \mathbb{R}^{n}$ is the 
system's state (random process), $\bm u(t)\in \mathbb{R}^{m}$ is the 
(deterministic) control, and $\bm x_0(\omega)\in \mathbb{R}^n$ 
is an random initial condition 
with prescribed probability density function $p_0(\bm x)$. As is 
well known \cite{Venturi_prs,venturi_mz,venturibook}, 
uncertain parameters in $\bm f$ can always be transferred to 
the initial condition. The solution to the Cauchy problem 
\eqref{eq:dynamics} is a function of the initial 
state $\bm x_0$ and a functional of the 
control $\bm u(t)$, i.e., 
$\bm x(t)=\bm x(t,\bm x_0, [\bm u])$.
We aim at designing $\bm u(t)$ by solving the optimal control problem 
\begin{equation}
\begin{cases}
\displaystyle \min_{\bm u(t)} J([\bm x(t)],[\bm u(t)])
\vspace{0.1cm}\\ 
\quad \text{subject to:}\\
\begin{array}{ll}
\qquad \dot{\bm x}=\bm f(\bm x,\bm u), & \bm x(0)=\bm x_0(\omega) \quad \text{($\bm x_0$ random)}    \\
\qquad\bm g(\bm u(t))\leq \bm 0,  &   \textrm{(constraints on the control)}   \\
\qquad\bm h\left(\mathbb{E}\{\bm x(t)\}\right)\leq \bm  0, & \textrm{(ensemble path constraints)}
\end{array}
\end{cases}
\label{eq:controlPr}
\end{equation}
where  $J([\bm x(t)],[\bm u(t)])$ is a 
cost functional \cite{venturi2016numerical} 
that involves an expectation $\mathbb{E}\{\cdot\}$ over the 
probability distribution of the 
initial state $\bm x_0(\omega)$. 
For example, $J$ could be of the form
\begin{equation}
J([\bm x(t)],[\bm u(t)])=\mathbb{E} \{F(\bm x(t_f))\}+
\int_{0}^{t_f} r(\bm u(\tau))d\tau,
\end{equation}
where  $F(\bm x)$ and $r(\bm u)$ are smooth functions. 
Other examples of $J$ involve 
the probability of detecting an unknown target
\cite{Automatica_search},  measures of 
risk \cite{Shapiro}, or indicator functions of 
disease propagation in a network of interacting 
individuals \cite{Keeling}. The optimal control 
problem \eqref{eq:controlPr} is, in general, a 
high-dimensional (large-scale) {\em non-convex} 
optimization problem. 
State-of-the-art algorithms to solve such problem are 
largely based on sampling the random initial state 
$\bm x_0$, e.g., using Gauss quadrature 
\cite{Automatica_search,JGCD_RS2015} or Monte 
Carlo methods \cite{Dick,SICON_Chris,venturi_ijhmt},  
and then computing a large ensemble of trajectories 
corresponding to a particular control $\bm u(t)$.
With such ensemble of paths available, one can then 
evaluate the performance metric $J$ by approximating the 
expectation operator with an appropriate cubature rule 
\cite{Automatica_search,SICON_Chris,walton_IFAC2014,
walton_NOLCOS2016}, and close the optimal control loop
with a suitable optimizer (see Figure \ref{fig:sketch}).
The formulation of optimal control strategies 
based on sample trajectories is straightforward, 
and it yields an optimal control problem 
which can be solved by classical algorithms, e.g., 
pseudospectral methods \cite{CMP_JGCD01,PS-TAC,PS-JOTA}.  
However, such algorithms are computationally intensive, 
and they are not effective even in a moderate number of 
state variables for the following reasons:  First, the set  
of sample trajectories can propagate in phase space 
in a very complicated manner. This implies that a large 
number of paths is usually required to evaluate the performance 
metric. Secondly, as the number of sample paths increases, 
the dimension of the discretized optimal control 
problem rapidly becomes intractable.
In particular, for a dynamical system with $n$ state variables 
and $M$ sample paths, the nonlinear control problem has dimension 
$n M$ ($M$ copies of the $n$-dimensional dynamical system
at each iteration of the optimizer). 
This is one of the main bottlenecks that 
limits the applicability of probabilistic 
collocation methods to systems with 
low dimensional random input vectors. For systems with 
large number of random variables,  
the memory requirements and computational cost 
becomes unacceptably large.

In this paper, we aim at overcoming these 
limitations by developing a  {\em new 
algorithm for computational optimal 
control under uncertainty} that is 
applicable to high-dimensional stochastic 
dynamical systems, and engineering control 
applications requiring online (real-time) implementations. 
The new algorithm is based on interior point optimization 
methods \cite{Biegler,biegler2009large}, common sub-expression elimination, 
and exact gradients obtained via automatic differentiation and 
computational graphs \cite{abadi2016tensorflow,baydin2018automatic}. 
They key features of the proposed algorithm are: 

\begin{enumerate}
\item {\em  Multi-shooting optimal control schemes}:
Multi-shooting is known for its numerical stability in 
solving deterministic optimal control problems. Based on 
piece-wise constant (in time) control approximation, we modified the  
deterministic algorithm and made it effective for uncertain 
optimal control problems. This includes imposing particle-independent 
continuity conditions  across adjacent time segments.
		
\item {\em Accurate gradient computations}: The optimization 
solver we developed combines an interior point method 
\cite{Biegler,biegler2009large,xu2014pyipopt} which converges most efficiently 
using accurate gradient calculations over the controlled 
ensemble.  In this work we compare a number of established 
techniques for gradient computation, namely traditional 
operator overloading algorithmic differentiation (ADOLC) 
\cite{walther2003adol}, more modern graph based 
methods for gradient calculation (TensorFlow) 
\cite{abadi2016tensorflow} and the most commonly 
used finite differences approach.  We demonstrate that 
graph based methods offer precision with considerable 
speed improvement over more traditional algorithms, 
albeit at the expense of considerable memory utilization, 
and slowness in initialization.

\item {\em Efficient memory utilization}:  
The algorithm has extremely low memory 
requirements, which means that it allows us to 
process a massive number of sample trajectories 
(in parallel) and determine the performance metric 
and associated optimal controls in a very efficient way. 
To achieve these results we developed a common 
sub-expression elimination (CSE) 
technique that can substantially reduce memory 
consumption during computations. As we shall see in 
section \ref{sec:CSE}, CSE considerably reduces the growth 
of memory requirements  with respect to the increase of 
the time discretization points in multi-shooting. 
This feature is {\em essential} for solving high dimensional 
problems and problems with long time horizon. 
It also reduces the discretization error in the 
control approximation, since smaller time steps can be afforded.  
Most importantly it does so at little to no loss to the performance 
of the gradient computation, and without any upfront initialization 
costs compared to graph based methods. 
\end{enumerate}

This paper is organized as follows. In section 
\ref{sec:multi-shooting} we formulate the optimal 
control problem \eqref{eq:controlPr} in a fully-discrete 
setting by using multi-shooting methods. This yields 
a non-convex optimization problem with nonlinear constraints 
which can be solved using interior point methods. To this end, 
it is extremely useful to have available fast and accurate 
gradient information, 
which we address in section \ref{sec:TensorFlow-IPOPT} 
using dataflow graphs, automatic differentiation, and 
common sub-expression elimination. 
In section \ref{sec:V&V} we develop a new 
criterion for verification and validation of the computed 
optimal controls using Pontryagin's minimum principles. 
In section \ref{sec:numerics} we demonstrate the 
accuracy and computational efficiency of the proposed 
control algorithms in applications to stochastic path-planning 
problems involving Unmanned Ground Vehicles (UGVs), 
Unmanned Aerial Vehicles (UAVs), and nonlinear PDEs. 
The main findings  are summarized in section \ref{sec:summary}.

\section{Multi-shooting schemes for ensemble optimal control}
\label{sec:multi-shooting}
\begin{figure}[t!]
		\centering
		\includegraphics[width=13cm]{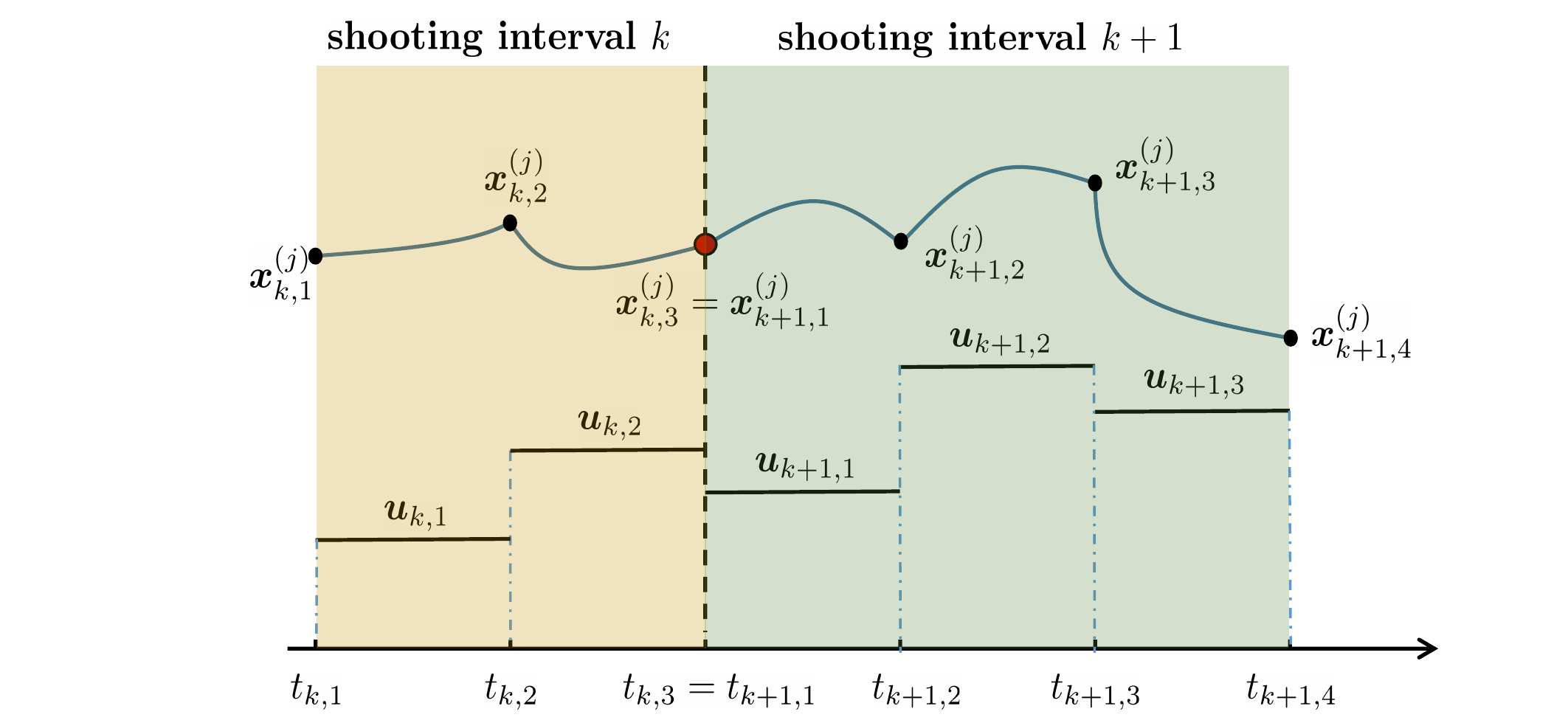}
\caption{Sketch of the multi-shooting setting for the optimal 
control of one-dimensional deterministic problem. Shown are the 
piece-wise constant approximation of the control $\bm u(t)$, and the 
piece-wise differentiable approximation of the state $\bm x(t)$.}
\label{fig:directshooting}
\end{figure}

\begin{figure}[t!]
\centering
	\includegraphics[width=13cm]{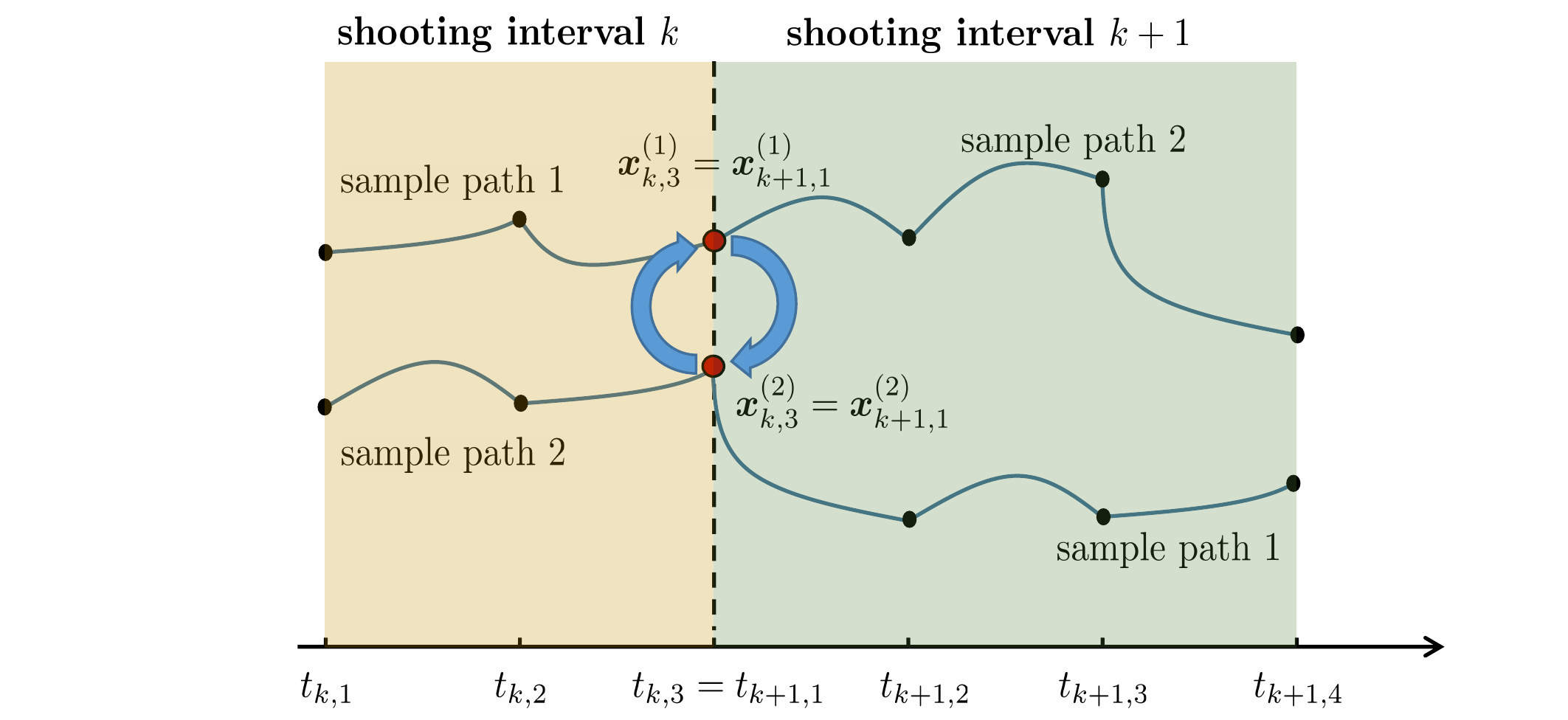}
\caption{Continuity of sample paths across shooting intervals. The continuity constraint 
\eqref{constr1} allows us to re-label sample paths at the boundary 
of each shooting segment. This does not affect ensemble statistical 
properties at all, but it massively reduces the number of 
path continuity constraints (equality constraints in the optimizer) 
from  $N$ to essentially one.}
\label{fig:trajectory_labels}
\end{figure}

Multi-shooting optimal control algorithms have been 
used extensively in control of deterministic systems 
because of their numerical stability and strightforward 
implementation.  In this section we describe an extension 
of multi-shooting to ensemble optimal control, i.e., control 
of ensembles of paths generated by the random dynamical 
system \eqref{eq:dynamics}. 
The first step when discretizing \eqref{eq:controlPr} with a
multi-shooting scheme is to divide control time 
horizon $[0,t_f]$ into $S$ sub-intervals, which 
we will call {\em shooting intervals}, 
\[[t_{k},t_{k+1}],\qquad  k=1,\ldots, S,\] 
with $t_1=0$ and $t_{S+1}=t_f$. 
Within each shooting interval we introduce an 
evenly-spaced grid of points
\begin{align*}
t_{k,j} =  t_{k}+(j-1)\Delta t_k,\qquad j=1,2,\ldots, N_k,  
\quad  k=1,\ldots, S,
\end{align*}
where $\Delta t_k$ denotes the time step size 
used in the $k$th shooting interval $[t_{k},t_{k+1}]$. 
Note that different shooting intervals can have different 
step size $\Delta t_k$, depending on the size of the 
temporal grid, i.e., $N_k$. The total number 
of time discretization points is 
\begin{align}
N=\sum_{k=1}^S N_k.
\label{Nt}
\end{align}
Within each shooting interval $[t_{k},t_{k+1}]$, 
the control $\bm u(t)$ is approximated by 
a {\em piecewise constant function} over 
the grid $\{t_{k,1},\dots, t_{k,N_k}\}$.
At this point, it is convenient to introduce the following notation:
\begin{equation}
\label{eqn:notation}
\bm x^{(i)}_{k,j} = \bm x^{(i)}(t_{k,j})\qquad \bm u_{k,j}=\bm u(t_{k,j})
\end{equation}
for the discretized trajectory of the state vector and control,
where $i=1,\ldots, M$ labels a specific realization 
of state process ($M$ sample paths total),
$k=1,\ldots,S$, labels the shooting interval 
$[t_{k},t_{k+1}]$, and $j=1,\ldots,N_k$ labels the 
time instant within the $k$th shooting interval. 
In Figure \ref{fig:directshooting}, we summarize the 
notation we used for the approximation/discretization of 
the state vector $\bm x(t)$ and the control $\bm u(t)$.
Note that piece-wise constant controls yield 
continuous sample paths $\bm x^{(i)}(t)$ 
with cusps at $t_{k,j}$.
With the piece-wise constant control approximation 
available, a sample of the state vector at any time instant 
can be computed by numerical integration as
\begin{equation}\label{eq:multishooting_x}
		\bm x^{(i)}_{k,j+1} \simeq \bm x^{(i)}_{k,j}+\int_{t_{k,j}}^{t_{k,{j+1}}} \bm f(\bm x^{(i)}(\tau), \bm u_{k,j}) d\tau, 
		\qquad i=1,\dots,M
\end{equation}
$M$ being the total number of sample paths.
%
%
%
%
To ensure the continuity of each  path across 
different shooting intervals, one can impose  
the continuity constraints
\begin{align}
\bm x^{(i)}_{k+1,1} =\bm x^{(i)}_{k,N_{k}}, \qquad i=1,\ldots,M.
\label{pointwise}
\end{align}
However, this formulation introduces a very 
large number of optimization constraints (one for each 
sample path), making the discrete optimization problem 
very hard to solve. In fact, a large number constraints 
usually increases the computational cost, and the size 
of unfeasible regions. Therefore, instead of adding one 
continuity constraint for each sample path, we introduce 
a single constraint for each dimension and each partition, 
but across all sample paths. 
\begin{equation}
\frac{1}{M} \sum_{l=1}^M 
\left(\bm x^{(l)}_{ k, N_{k}} -\bm x^{(l)}_{k+1, 1}\right)^2=0
\label{constr1}
\end{equation}
which still guarantees continuity of sample 
trajectories across different shooting intervals, 
but in a way that is  {\em agnostic about the 
trajectory labels}.
In other words, across different shooting intervals
the constraint \eqref{constr1} allows for a reshuffling of the 
trajectories labels. As is well know, this does not affect 
the computation of ensemble statistical properties (see Figure \ref{fig:trajectory_labels}).
%
Approximating the integral in \eqref{eq:multishooting_x}, e.g., 
with an explicit $s$-stage Adams-Bashforth \cite{hairer} 
method yields 
\begin{align}
\bm x^{(i)}_{k,j+1} = \bm x^{(i)}_{k,j}+
\Delta t_k\sum_{q=0}^{s-1}b_j \bm f\left(\bm x^{(i)}_{k,j-s},\bm u_{k,j-s}\right).
\label{const2}
\end{align} 
At this point it is convenient to denote the $i$th sample of the 
fully discrete state process in $[0,t_f]$ as 
\begin{align*}
\bm X^{(i)}=\left\{\bm x^{(i)}_{1,1}, \ldots \bm x^{(i)}_{1,N_1},\bm x^{(i)}_{2,1}, \ldots, \bm x^{(i)}_{2,N_2},\ldots, \bm x^{(i)}_{S,1}, \ldots, \bm x^{(i)}_{S,N_S}\right\}, \qquad i=1,...,M.
\end{align*}
It is also convenient to define a vector 
collecting the sample values the process $\bm x(t)$ 
at the boundaries of the shooting intervals, i.e., 
\begin{align*}
\bm X_b^{(i)}=\left\{\bm x^{(i)}_{1,1},\bm x^{(i)}_{1,N_1},\bm x^{(i)}_{2,1},\bm x^{(i)}_{2,N_2}, \ldots, \bm x^{(i)}_{S,1}, \bm x^{(i)}_{S,N_{S}}\right\}, 
\end{align*}
and the full vectors of states and boundary values
\begin{align*}
\bm X = \left\{\bm X^{(1)}, \ldots, \bm X^{(M)}\right\}, \qquad 
\bm X_b = \left\{\bm X_b^{(1)}, \ldots, \bm X_b^{(M)}\right\}.
\end{align*} 
Similarly, we denote the discrete control vector as 
\begin{align*}
\bm U=\left\{\bm u_{1,1}, \ldots \bm u_{1,N_1-1}, \bm u_{2,1},\ldots, \bm u_{2,N_{2}-1},\ldots,\bm u_{S,1},\ldots, \bm u_{S,N_{S}-1}\right\}.
\end{align*}
With this notation, we can write the fully discrete form of the 
optimal control problem \eqref{eq:controlPr} as 
\begin{equation}
\begin{cases}
\displaystyle \min_{\bm U,\bm X_b} J(\bm U,\bm X_b)\\ 
\quad \text{subject to:}\\
\displaystyle \qquad\qquad \bm x^{(i)}_{k,j+1} = \bm x^{(i)}_{k,j}+
\Delta t_k\sum_{q=0}^{s-1}b_j \bm f\left(\bm x^{(i)}_{k,j-s},\bm u_{k,j-s}\right),
\quad \text{$\bm x_{1,1}^{(i)}$ random}
\qquad\\
\displaystyle \qquad\qquad 
\sum_{k=1}^S \sum_{i=1}^M 
\left\|\bm x^{(i)}_{ k, N_{k}} -\bm x^{(i)}_{k+1, 1}\right\|_2^2=0
\qquad \text{(path continuity constraint)}
\vspace{0.2cm}\\
\displaystyle \qquad\qquad \bm g(\bm U)\leq \bm 0 
\qquad \textrm{(constraints on the control)}\\
\displaystyle \qquad\qquad  \bm h(\mathbb{E} \{ \bm X \}) \leq \bm 0
\qquad \textrm{(ensemble path constraints)}
\end{cases}
\label{eq:discrete_opt}
\end{equation}

{\remark $\,$} The first constraints in \eqref{eq:discrete_opt} 
is the discrete form of the ODE \eqref{eq:dynamics}, i.e., it defines 
the temporal dynamics of each sample path under the control $\bm U$. 
In practice, all paths are pushed forward in time 
{\em simultaneously} for a given control by vectorizing the state 
vector in a way that includes all solution samples. 
Marching forward in time can take advantage of the massively 
parallel structure of the discretized ODE system.

\section{Fast gradient computations}
\label{sec:TensorFlow-IPOPT}
The large-scale {constrained optimization problem} 
\eqref{eq:discrete_opt} can be solved by using 
optimization algorithms based on gradient information. 
To this end, it is very important to be able to compute 
the gradient of the discretized cost function $J(\bm U,\bm X_b)$ 
with respect to the decision variables $(\bm U,\bm X_b)$, 
with accuracy and efficiency. In this section we provide 
technical details on how we implemented such gradient 
computation using automatic differentiation functions 
available in graph-based methods such as TensorFlow 
\cite{abadi2016tensorflow}, and backward propagation. 
Our algorithm leverages the fact that forward time integration 
can be framed as a recurrent neural network, 
hence providing {exact gradients} at a negligible 
computational cost. Indeed, as we will see, all 
operations can be performed very efficiently on  
computational graphs.

\subsection{Data-flow graphs}  
\label{sec:dataflowGraphs}
\begin{figure}[!t]
	\begin{subfigure}[b]{.48\linewidth}
		\centering
		\includegraphics[width=\linewidth]{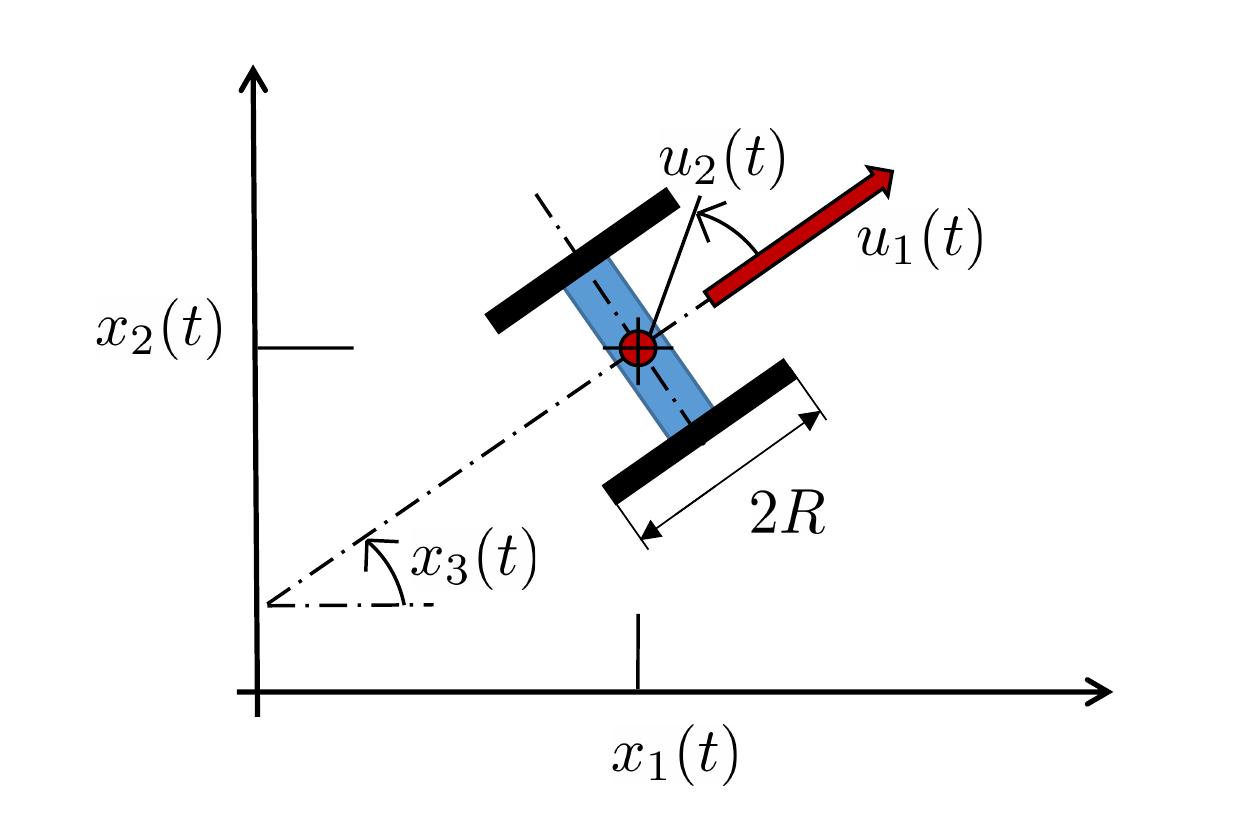}
		\caption{\footnotesize UGV modeled by the system \eqref{example:UGV}.}
		\label{fig:ugv_opt_unc}
	\end{subfigure}
	\begin{subfigure}[b]{.48\linewidth}
		\centering
		\includegraphics[width=\linewidth]{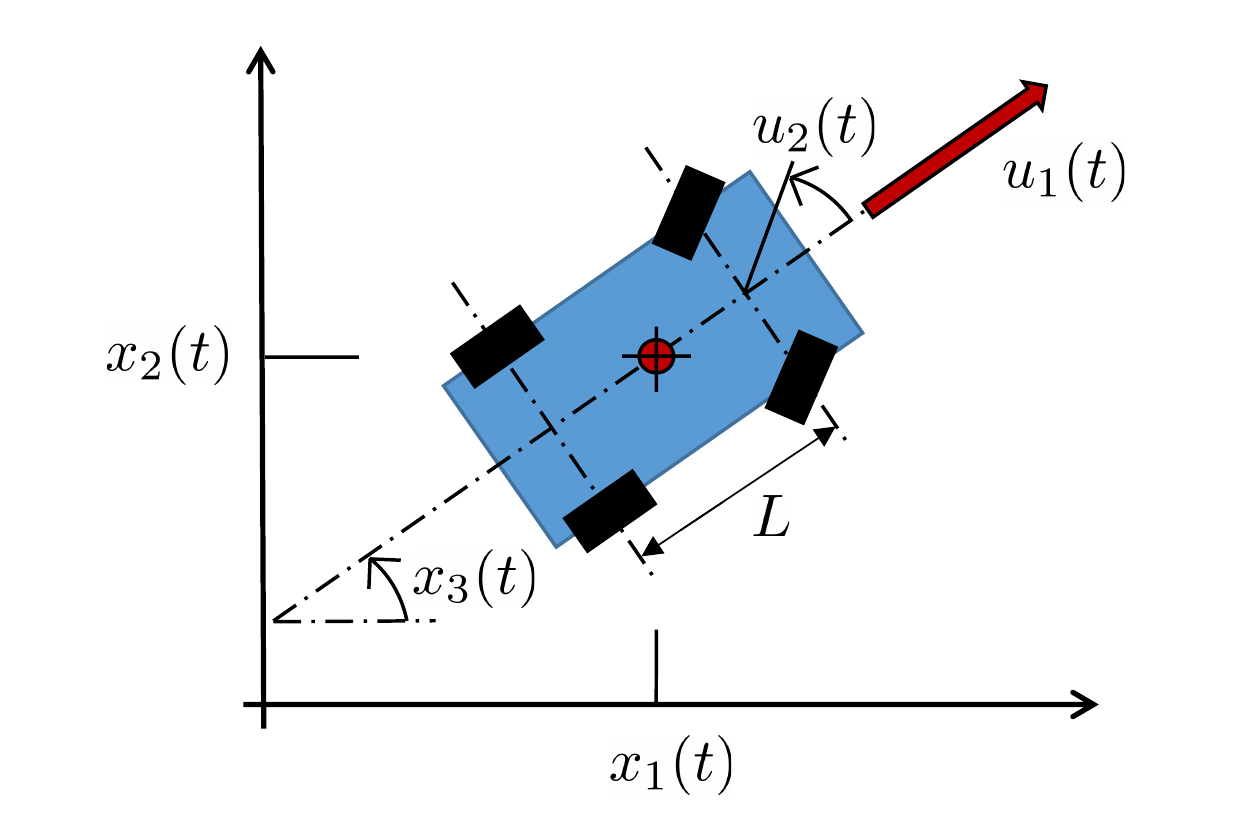}
		\caption{\footnotesize  UGV modeled by the system \eqref{eq:UGV2}.}
		\label{fig:ugv_nom_unc}
	\end{subfigure}
	\caption{Sketches of two different Unmanned Ground Vehicles (UGVs). Shown are the state space variables $\{x_i(t)\}_{i=1,2,3}$, and the controls $u_1(t)$ and $u_2(t)$. }
	 \label{fig:rear_car}
\end{figure}

\begin{figure}[!t]
	\centering 
	\includegraphics[width=11cm]{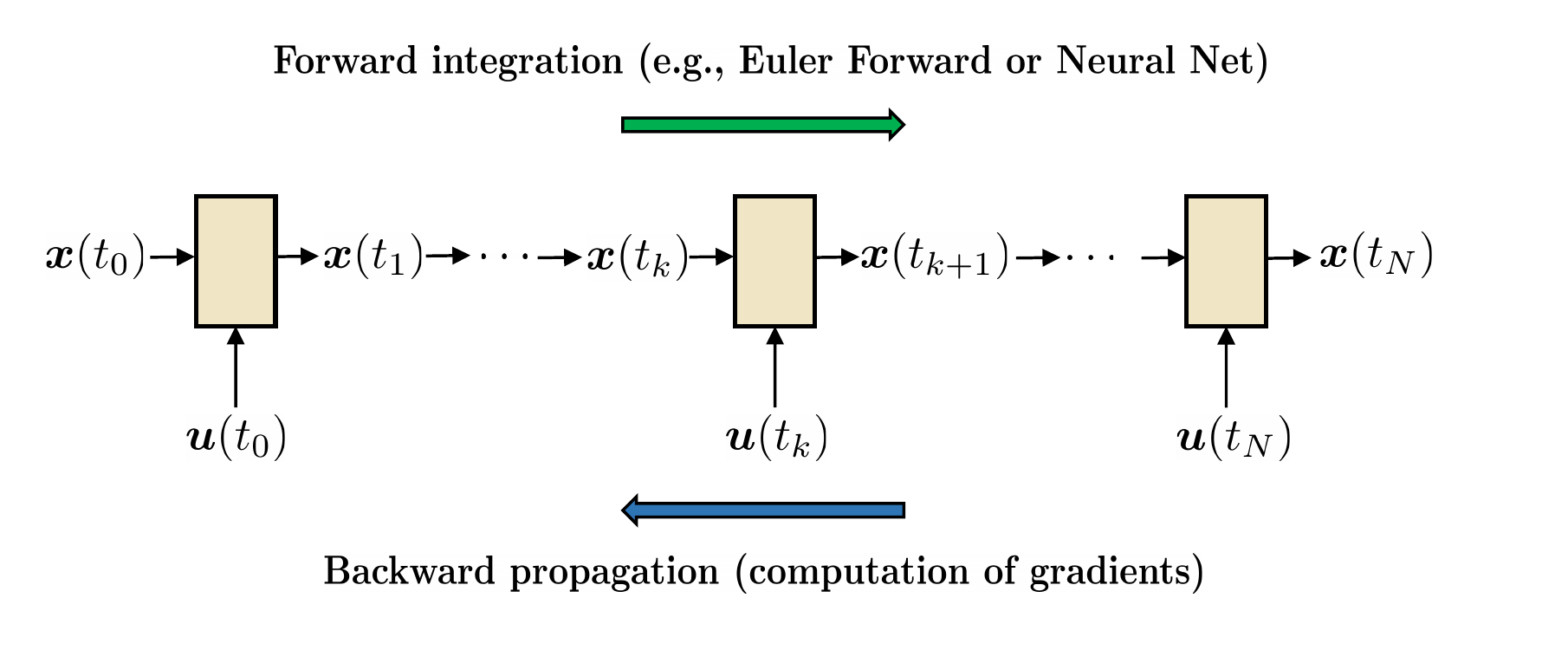}\\
	\includegraphics[width=18cm]{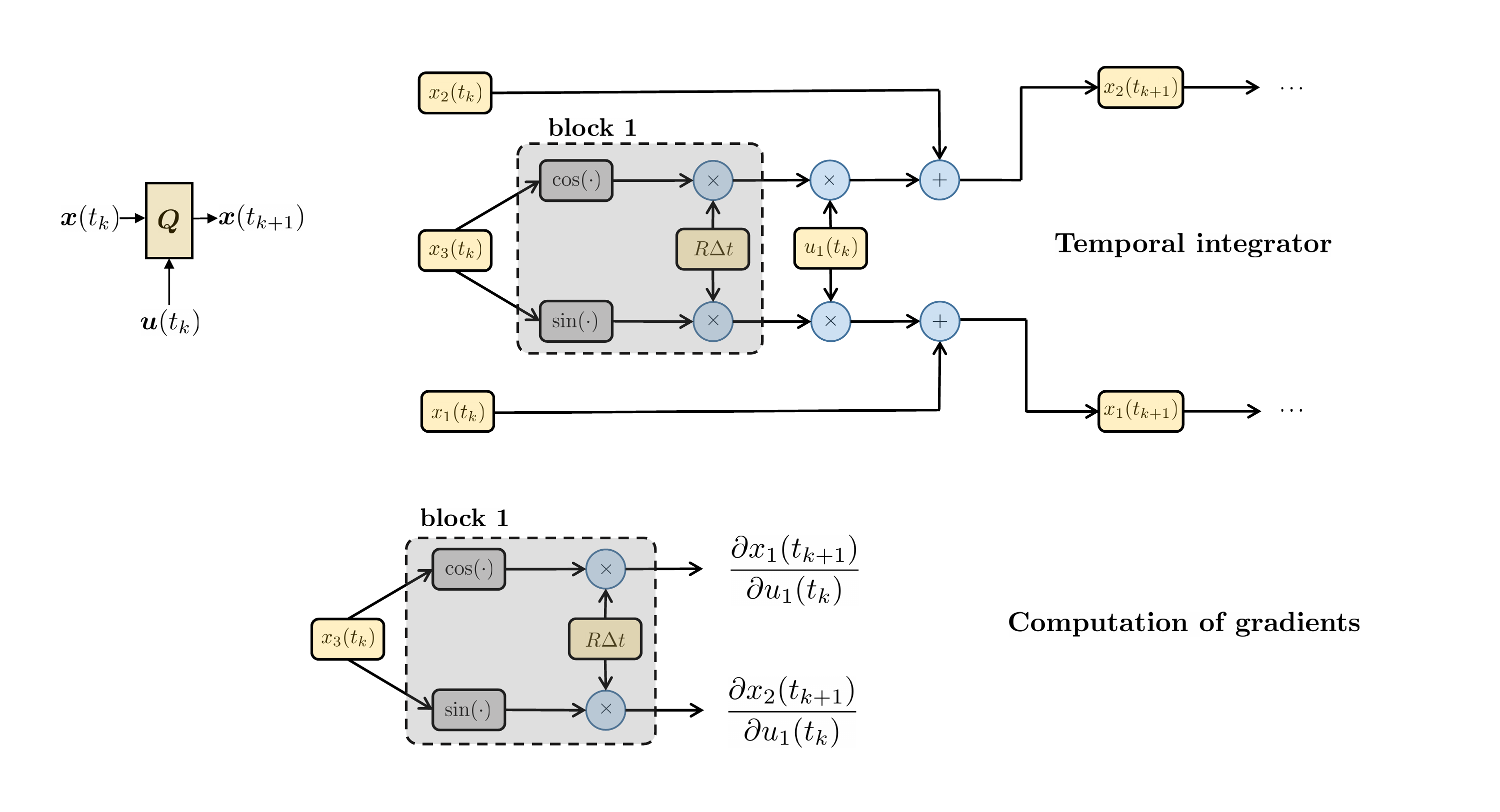}
	\caption{Portion of the data graph for the UGV model discretized with 
		Euler forward time integration and graph the compute the 
		kth component of the gradient of the cost functional. Note that 
		there are sub-graphs that are shared among the computations 
		(boxed graphs). Using {reactive programming}, all operations 
		within such boxes are computed only { once}, and then shared 
		across the graph.}
	\label{fig:ugvEuler}
\end{figure}

Modern machine learning techniques use 
decentralized data graphs to map computations in
different nodes of a computer cluster 
\cite{abadi2016tensorflow}.  Edges and nodes of such graph usually 
represent flows of data and input-output maps (operations), 
respectively. The dataflow graph is pre-compiled, optimized and 
stored as metadata in memory for the duration of the calculation.  
Using reverse automatic differentiation (backpropagation), the 
computational graph of the gradients with respect to a given cost 
function is simultaneously constructed and similarly stored.  
Using these two graphs one can integrate the system 
in time for a given control and compute 
gradients of the state trajectories with respect to the 
controls. The calculations are performed using a 
data event driven mechanism allowing for nearly 
simultaneous calculations of gradients as the dynamical 
system is being integrated.  
To illustrate the main idea, consider the follow 
dynamical system modeling the wheel-driving 
Unmanned Ground Vehicle (UGV) 
shown in Fig.\ref{fig:rear_car}(a).
\begin{equation}\label{example:UGV}
\begin{dcases*}
\dot{x}_1= Ru_1(t) \cos(x_3) \\
\dot{x}_2= Ru_1(t) \sin(x_3) \\
\dot{x}_3= R u_2 
\end{dcases*}
\end{equation}
Here, $(x_1(t),x_2(t))$ denotes the position of the vehicle 
on the Cartesian plane while $x_3(t)$ is the  heading angle. 
The controls are $u_1(t)$ (velocity) and $u_2(t)$ (steering angle). 
For illustration purposes, we discretize \eqref{example:UGV} 
with the Euler forward scheme (one-step Adams Bashforth 
\eqref{const2}). This yields,
\begin{equation}
\label{example:UGVscheme}
\begin{dcases*}
x_1(t_{k+1})=x_1(t_k)+ \Delta t R u_1(t_k) \cos(x_3(t_k)) \\
x_2(t_{k+1})=x_2(t_k)+ \Delta t R u_1(t_k) \sin(x_3(t_k)) \\
x_3(t_{k+1})=x_3(t_k)+ \Delta t R u_2(t_k) 
\end{dcases*}.
\end{equation}
This system of equations defines explicitly an input-output 
map of the form\footnote{Explicit Runge-Kutta and linear 
multistep methods can always be written in the 
form \eqref{eq:onestep} (see \cite{Bram2019,hairer,Reddy1992}).}
\begin{equation}
\bm x_{k+1} = \bm Q(\bm x_{k} ,\bm u_k).
\label{eq:onestep}
\end{equation}
Such map takes in 
$\bm x_k=[x_1(t_k),x_2(t_k),x_3(t_k)]$ and 
$\bm u_k=[u_1(t_k),u_2(t_k)]$ (which we assumed constant 
within the time interval $\Delta t$), and sends them to 
$\bm x(t_{k+1})$. Such nonlinear map can be conveniently 
represented as a ``black box'' with inputs 
$\bm x(t_k)$ and $\bm u(t_k)$, 
and output $\bm x(t_{k+1})$ as illustrated in Figure \ref{fig:ugvEuler}.
It is clear at this point that the process of integrating 
the state vector $\bm x(t_0)$ forward in time 
from $t_0$ to $t_N=t_f$ for any given piece-wise 
constant (in time) control
\begin{align*}
\bm U = \left[ \bm u(t_0),\cdots \bm u (t_{N-1})\right]
\end{align*}
can be seen as a composition of the same nonlinear map 
\eqref{example:UGVscheme}. This essentially defines a 
recurrent network that allows for an extremely low 
memory footprint and fast gradient computation.
Moreover, as clearly seen from the computational graph sketched 
in Figure \ref{fig:ugvEuler}, portions of the calculations for 
the forward trajectory propagation as well as the backward gradient 
calculation are shared.  Using {reactive programming} such 
shared portions of the graphs are computed only once, 
while {memoization} realizes additional substantial 
computational savings. 
The exact amount of savings depends on the efficiency of the graph 
compiler as well as the particular dynamics. Once the compiler 
creates and compiles the data graph in memory for a single 
trajectory, adding a dimension for sampling does not involve 
any changes to the graph or gradient calculations.  By choosing 
a non-adaptive integration scheme, we ensure that the operations 
can be performed in lock step mode across all samples, and can 
utilize hardware with wide SIMD capabilities such as 
GPU, TPU or ASIC cards.
%

\subsection{Memory reduction and accelerated computations via Common Sub-expression Elimination (CSE)}
\label{sec:CSE}
Graph based methods for gradient calculations utilized by Deep Neural Net frameworks are scalable and are implemented in a computationally efficient way.  However because they are designed to operate on large system for prolonged periods of time, ranging from hours to months, they have the following disadvantages:
\begin{enumerate}
\item They consume an inordinate amount of memory of the order of tens of gigabytes, especially for systems that have many constraints.

\item They use lazy initialization which makes the initial ramp up time very slow, often lasting hours for high dimensional systems with long time horizons.  This is acceptable for deep neural net architectures, but hampers the response time for solving high dimensional optimal control problem. 

\item Because of the hefty memory and time initialization requirements, they are unsuitable for online (real-time) applications, in particular those targeting autonomous vehicles.
\end{enumerate}
The main idea around common sub-expression elimination is 
the observation that during the construction of the computational 
graph in systems like Tensorflow, the memory increases linearly 
with the number of steps.  However since the computations 
are identical from on time step to the next, this motivates the 
objective of eliminating the common repeated calculations for 
each trajectory in order to save memory while not sacrificing 
computational speed, hence the name
common sub-expression elimination.

To illustrate the basic idea of CSE in a simple way, let us 
consider \eqref{eq:dynamics} and write the time-discrete 
form as a symbolic nonlinear one-step recursion of 
the form \eqref{eq:onestep}.
Suppose we are interested in minimizing a cost 
function of the form 
\begin{equation}
J(\bm U) = \frac{1}{M}\sum_{i=1}^M F(\bm x^{(i)}(t_f)),
\label{JD}
\end{equation}
where $F$ is smooth. 
The gradient of \eqref{JD} with respect to 
$\bm u(t_k)$ can be computed using  the 
chain rule as
\begin{equation}
\label{eqn:back_prop}
\frac{\partial J(\bm U)}{\partial \bm u(t_k)}= \frac{1}{M}\sum_{i=1}^M
\frac{\partial F(\bm x^{(i)}_{N})}{\partial \bm x}
\frac{\partial \bm Q(\bm x^{(i)}_{N-1},\bm u_{N-1}) }{\partial \bm x} 
\cdots 
\frac{\partial \bm Q(\bm x^{(i)}_{k+1},\bm u_{k+1}) }{\partial \bm x} 
\frac{\partial \bm Q(\bm x^{(i)}_{k},\bm u_k) }{\partial \bm u}. 
\end{equation}
This is the well-know {\em backward propagation formula} in 
recurrent neural networks. Indeed, the gradient of $J$ is 
obtained by iteratively applying the Jacobian maps 
$\partial \bm x (t_{j})/\partial \bm x (t_{j-1})$ from $j=N$ to $j=k+1$. Contrary to commonly used patterns in graph methods, here we 
need to store in memory only one gradient and two Jacobian functions, i.e., 
\begin{equation*}
\frac{\partial F}{\partial \bm x}, \qquad 
\frac{\partial \bm Q(\bm x,\bm u) }{\partial \bm x},
\qquad \text{and}\qquad 
\frac{\partial \bm Q(\bm x,\bm u) }{\partial \bm u},
\end{equation*}
where $\bm Q$ is defined in \eqref{eq:onestep}.
Even though other memory efficient schemes exist 
where the Jacobians are not stored explicitly, but rather the 
Jacobian vector product is done explicitly at each forward step, 
our numerical experiments indicate that 
storing the Jacobian explicitly strikes the best 
balance between computational efficiency and 
memory consumption.  Moreover, the Jacobians 
for dynamical systems are of much smaller size than the 
Jacobians for neural nets, especially when the number of 
hidden layers in the latter is large, which explains why Jacobians 
are not stored explicitly in neural net frameworks.
Taking advantage of the linear nature of the contribution 
of each trajectory of the ensemble to the end objective, we 
can leverage this linearity to customize the memory 
utilization through the mechanism of batching the samples.  
Specifically, a hundred-dimensional dynamical system over 
one thousand time steps requires a mere 100MB 
of memory per sample which is insignificant by 
today's standards.  On a modest embedded system 
or a GPU with 10GB or more in memory, this allows for 
the simultaneous propagation of one hundred 
trajectories per batch using SIMD instructions. 
For smaller systems, much larger batches of samples can 
be propagated simultaneously.  The specific number 
of optimal samples that can be propagated in a single 
batch needs to be determined by experimentation on 
the particular hardware architecture in question.  The 
overall number of samples propagated can therefore be 
as large as required by the problem in question, without 
increase in memory footprint.

\subsubsection{Common Sub-expression Elimination: An example}

\begin{figure}[t]
\centering
\includegraphics[height = 6.5 cm]{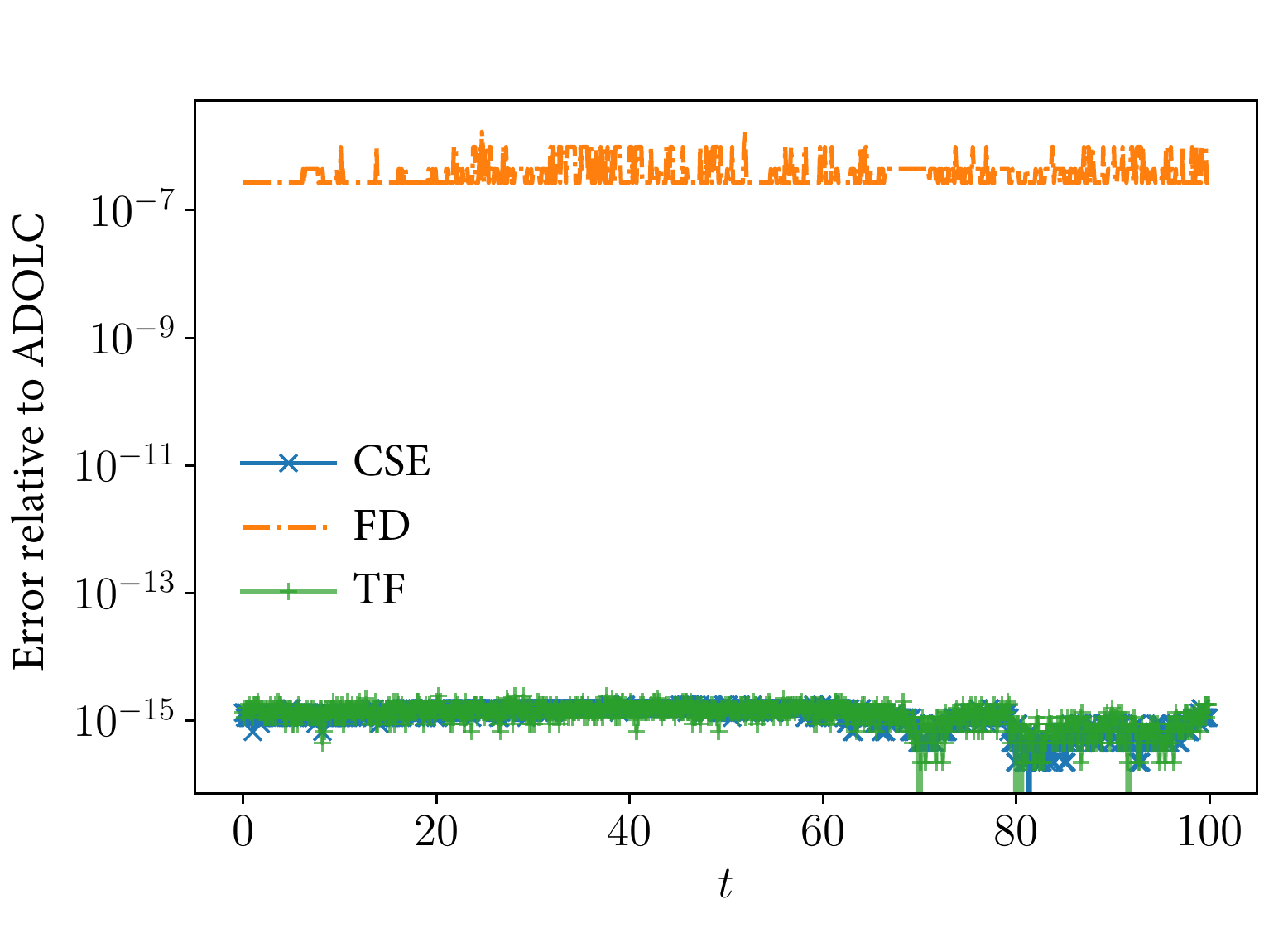}
\caption{Absolute error of $\partial J/\partial u_1(t)$ 
with respect to the benchmark gradient computed with ADOLC 
for CSE, second-order centered finite differences (FD) and TensorFlow (TF). It is seen that tensor flow and CSE are numerically exact.}
\label{fig:err1}
\end{figure}

\begin{figure}[t]
	\centering
	\begin{subfigure}[b]{.49\linewidth}
		\centering
		\includegraphics[height = 5.8cm]{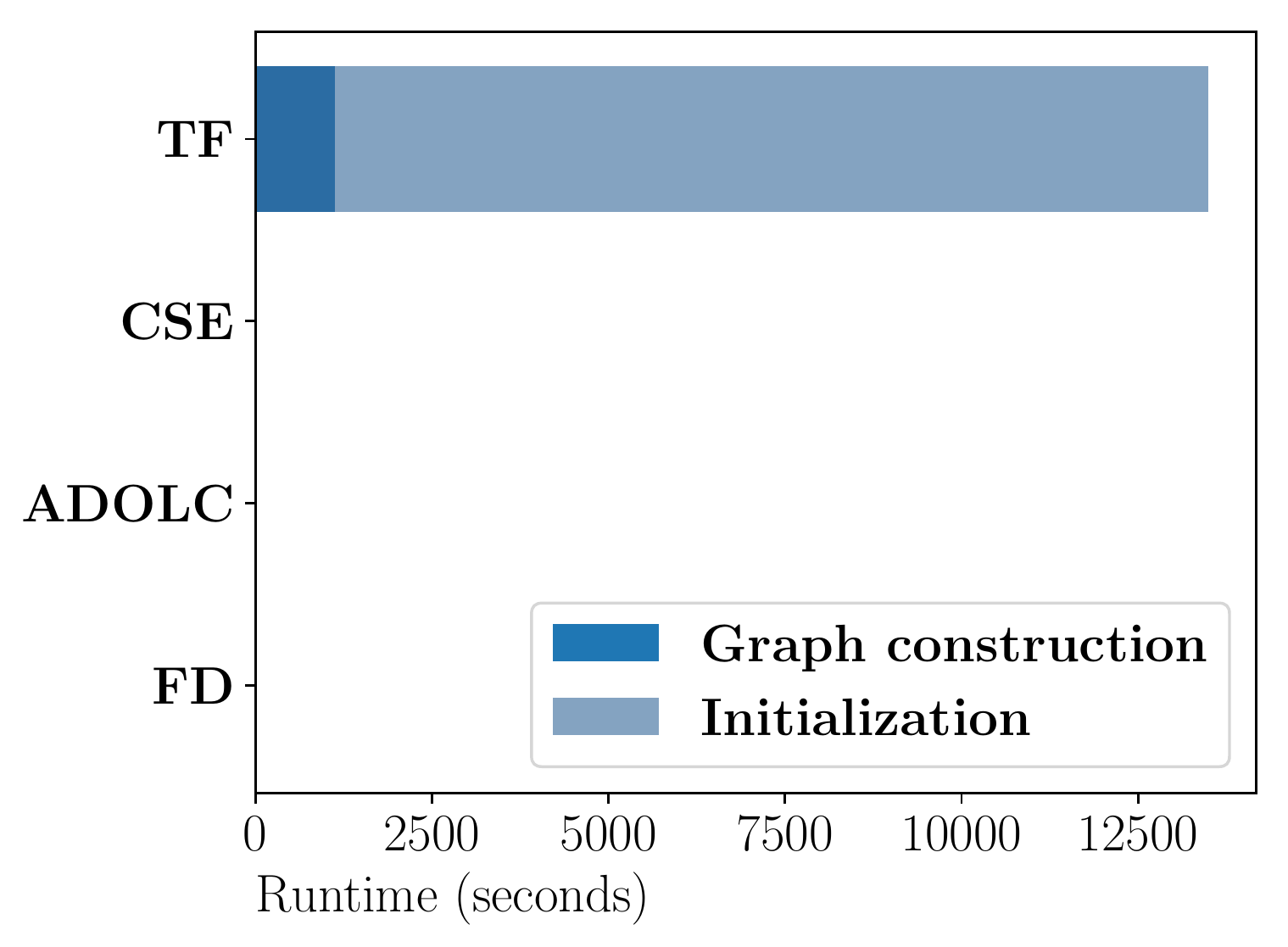}
	\end{subfigure}
	\begin{subfigure}[b]{.49\linewidth}
		\centering
		\includegraphics[height = 6.2cm,width=8.0cm]{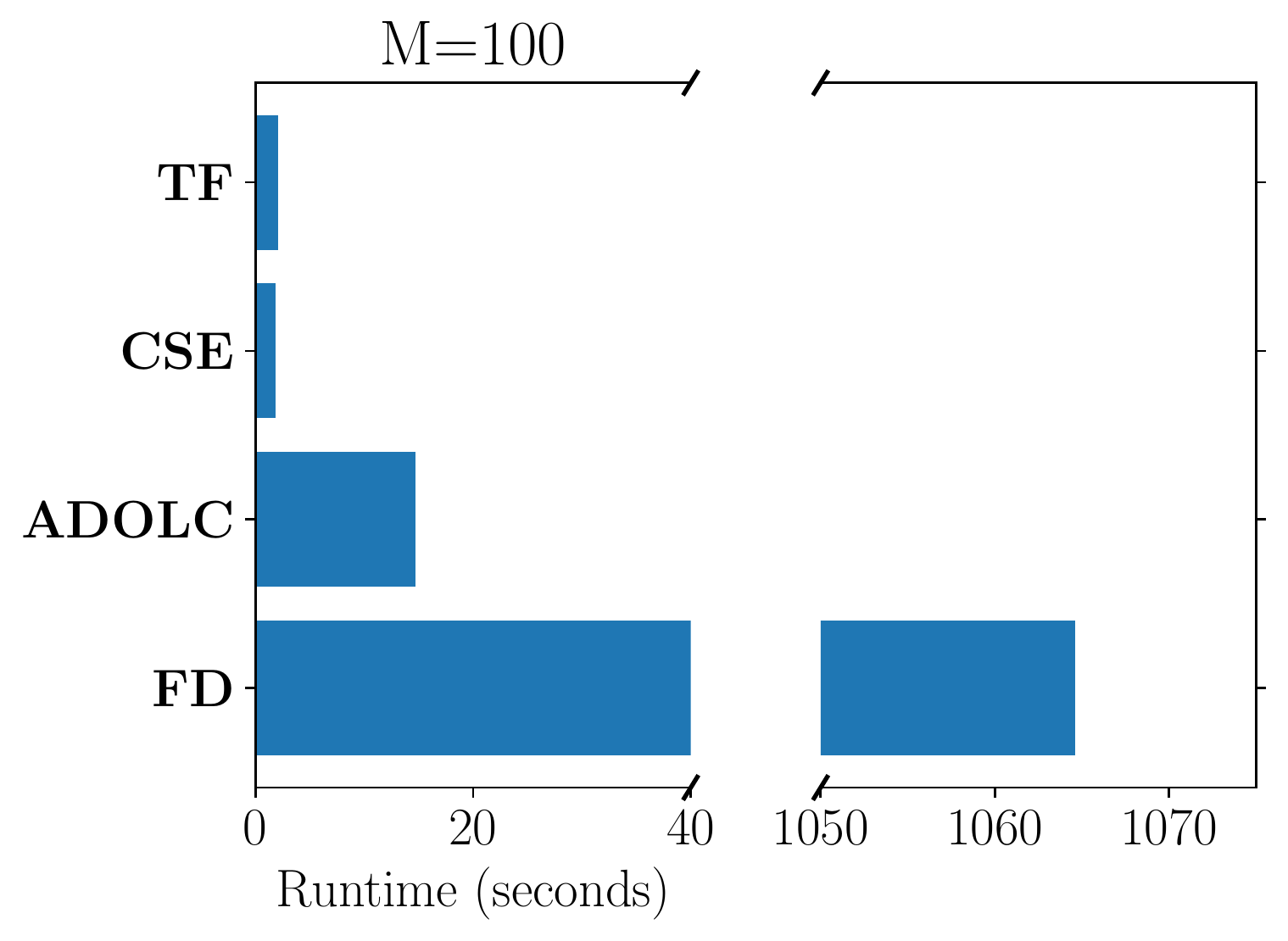}
	\end{subfigure}
	\begin{subfigure}[b]{.49\linewidth}
	\centering
	\includegraphics[height = 6 cm]{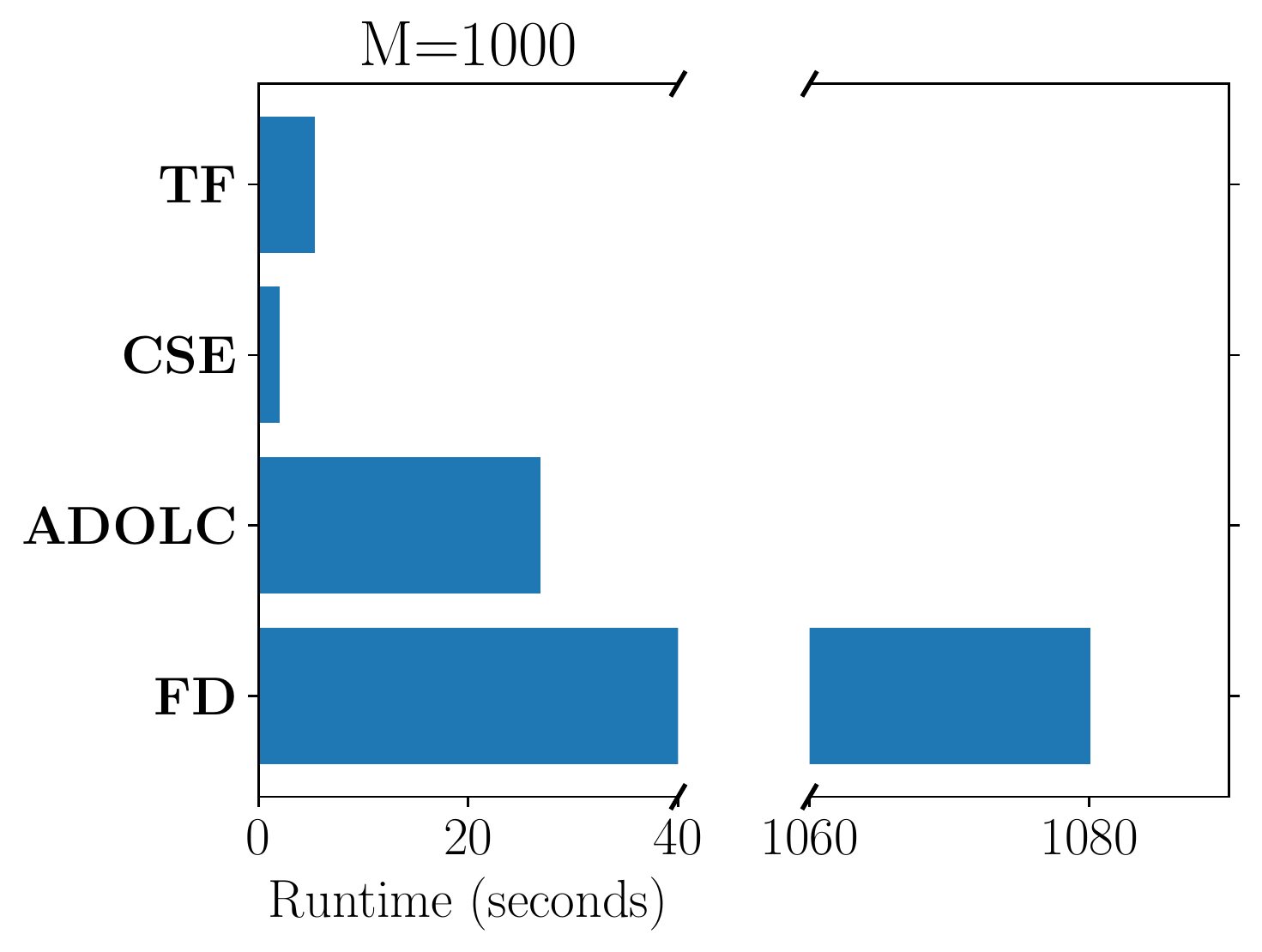}
	\end{subfigure}
	\begin{subfigure}[b]{.49\linewidth}
	\centering
	\includegraphics[height = 6 cm]{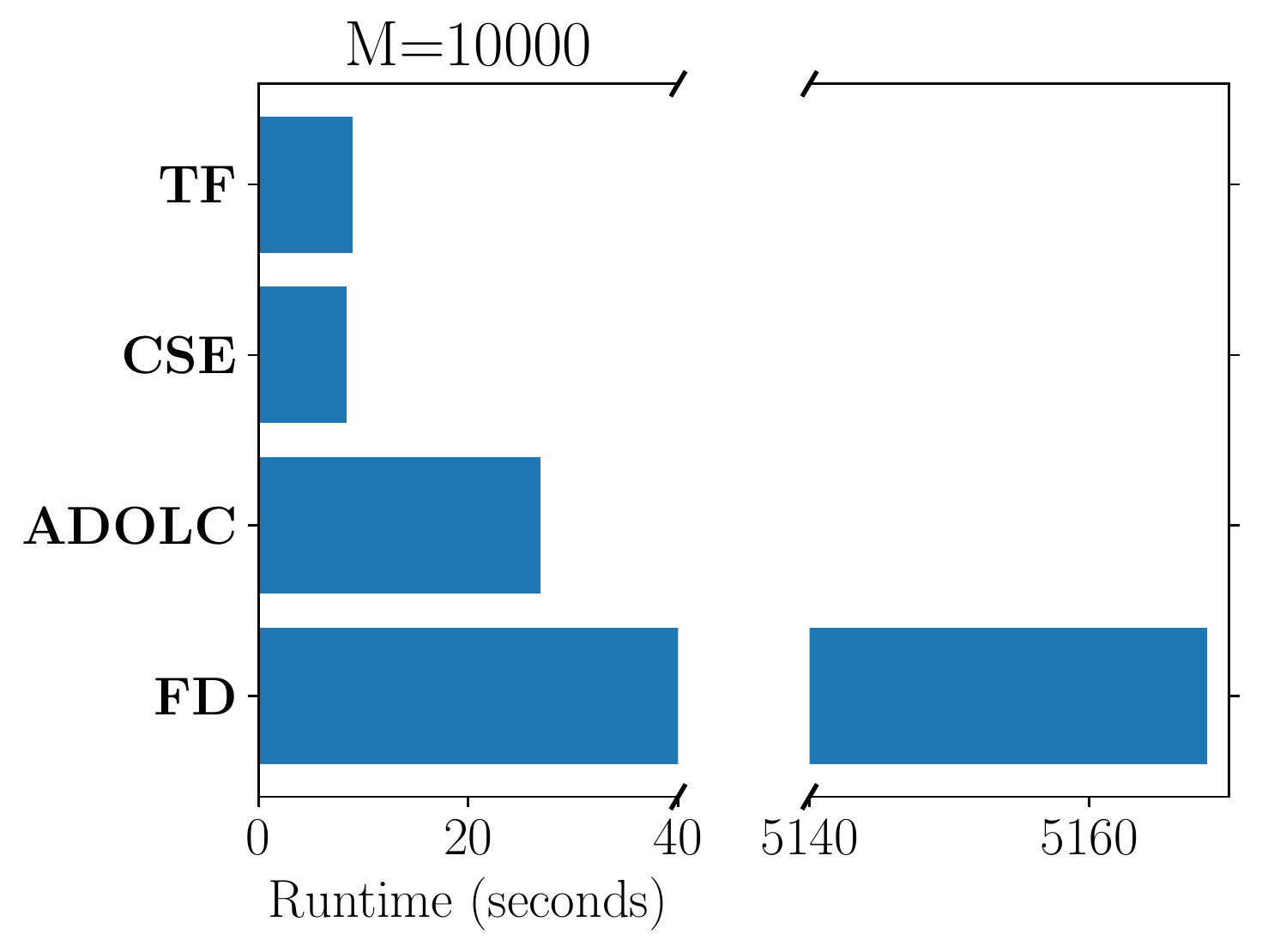}
	\end{subfigure}
\caption{Execution times in seconds 
for different gradient methods.  Here 
TF stands for tensorflow, and FD represents second-order finite 
differences.  The CSE algorithm uses batch sizes of 100 samples 
for each thread, whereas ADOLC is parallelized one sample at a time.  
Tensorflow automatically optimizes its configuration for the target 
hardware architecture.  Note the clear disadvantage in runtime 
that finite difference methods exhibit. All numerical experiments were conducted on a dual Xeon E5-2651 v2 processor workstation operating at 1.80GHz with 145 GB of memory.}
\label{fig:runtime}
\end{figure}

\begin{figure}[t]
	\centering
	\begin{subfigure}[b]{.49\linewidth}
		\centering
		\includegraphics[height = 6 cm]{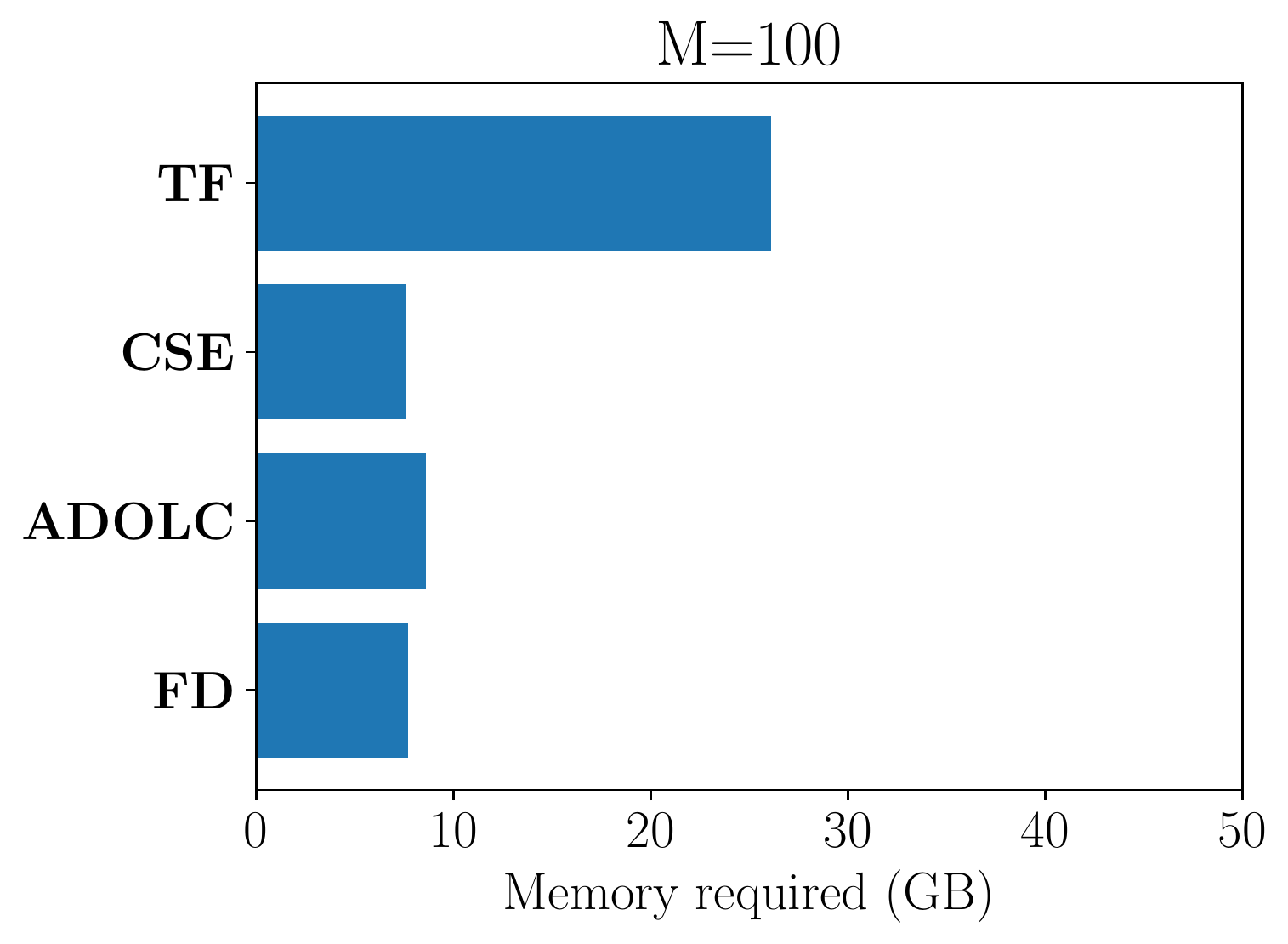}
	\end{subfigure}
	\begin{subfigure}[b]{.49\linewidth}
		\centering
		\includegraphics[height = 6 cm]{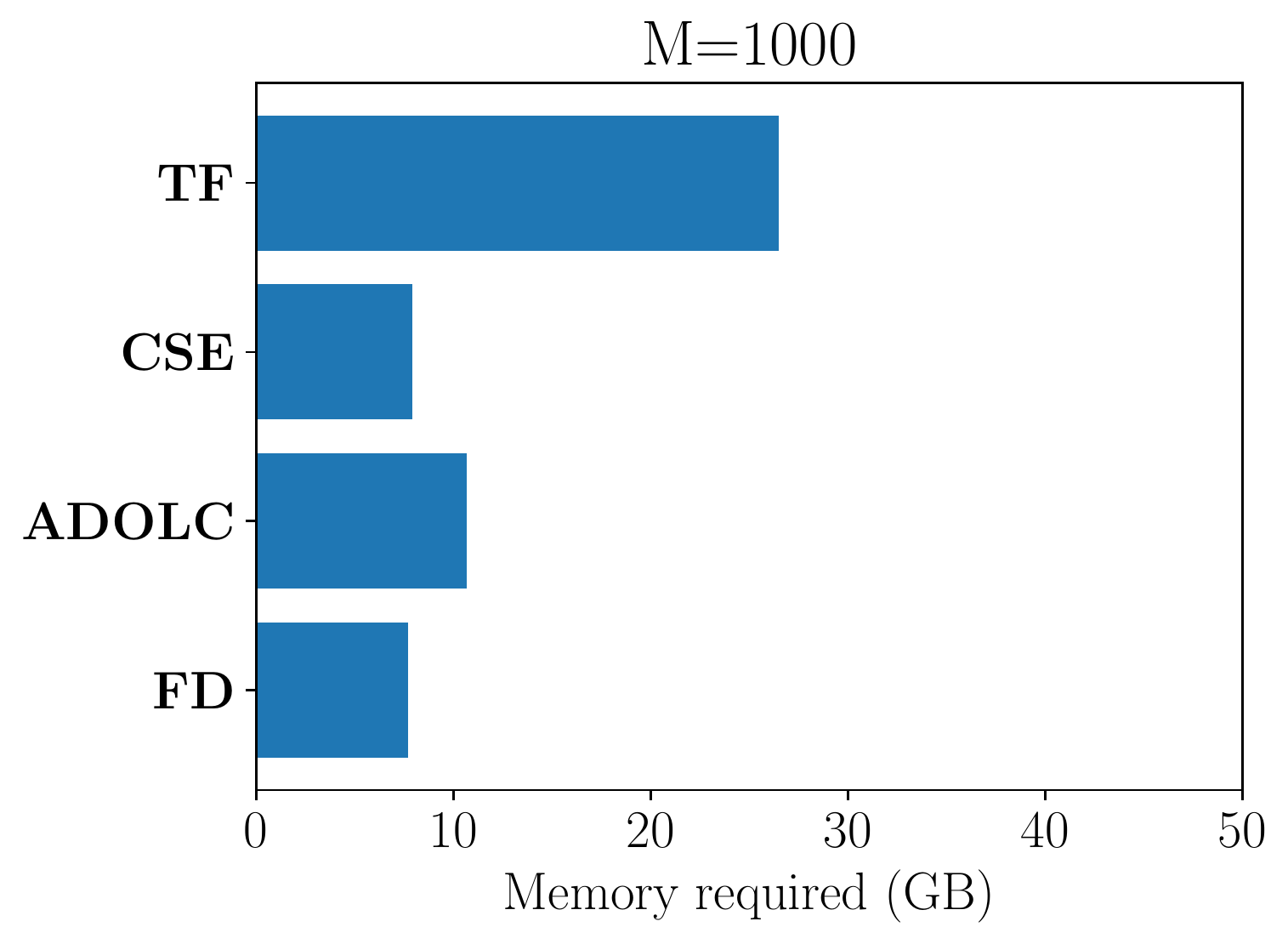}
	\end{subfigure}
	\begin{subfigure}[b]{.49\linewidth}
		\centering
		\includegraphics[height = 6 cm]{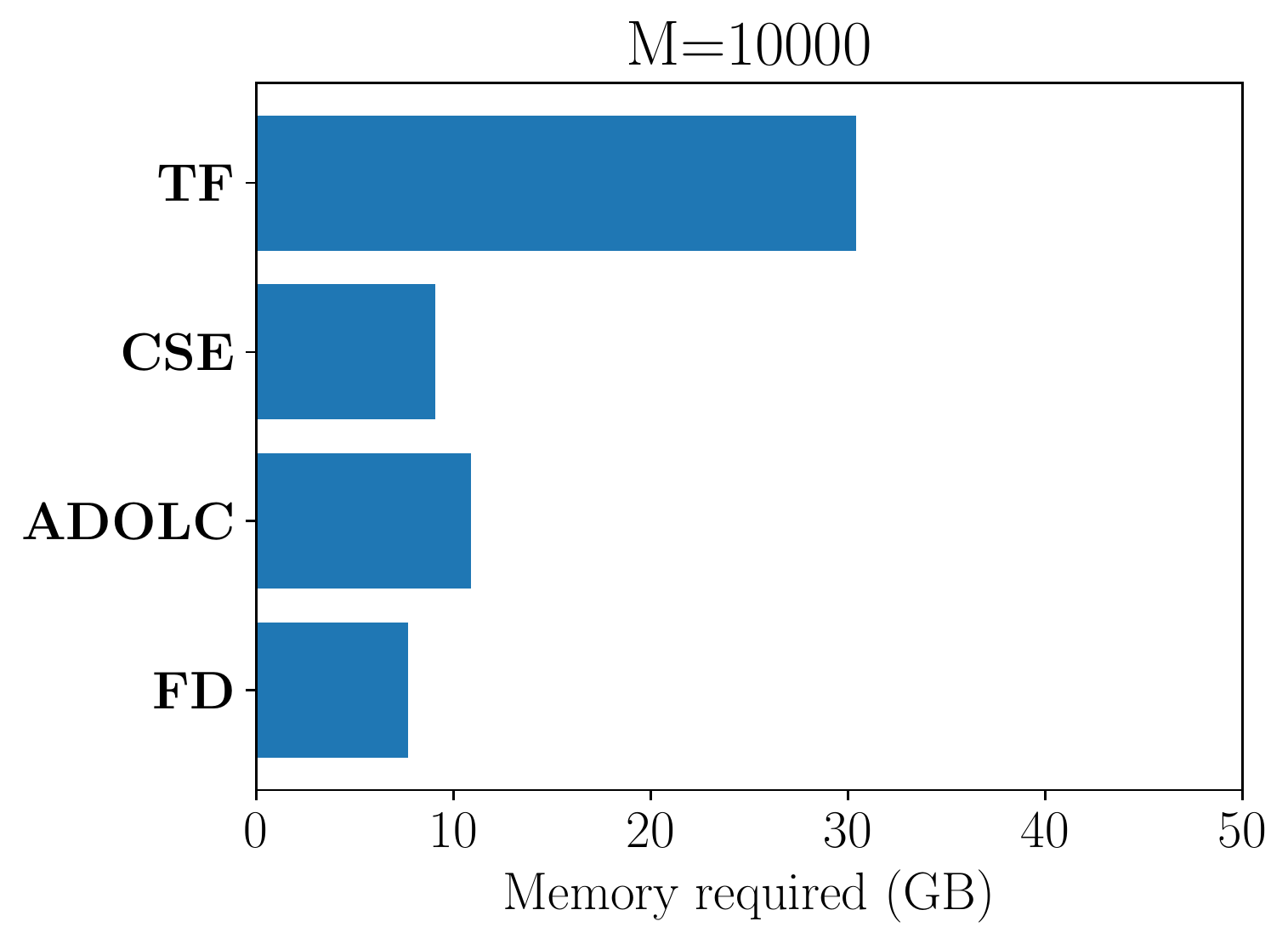}
	\end{subfigure}
	\begin{subfigure}[b]{.49\linewidth}
		\centering
		\includegraphics[height = 6 cm]{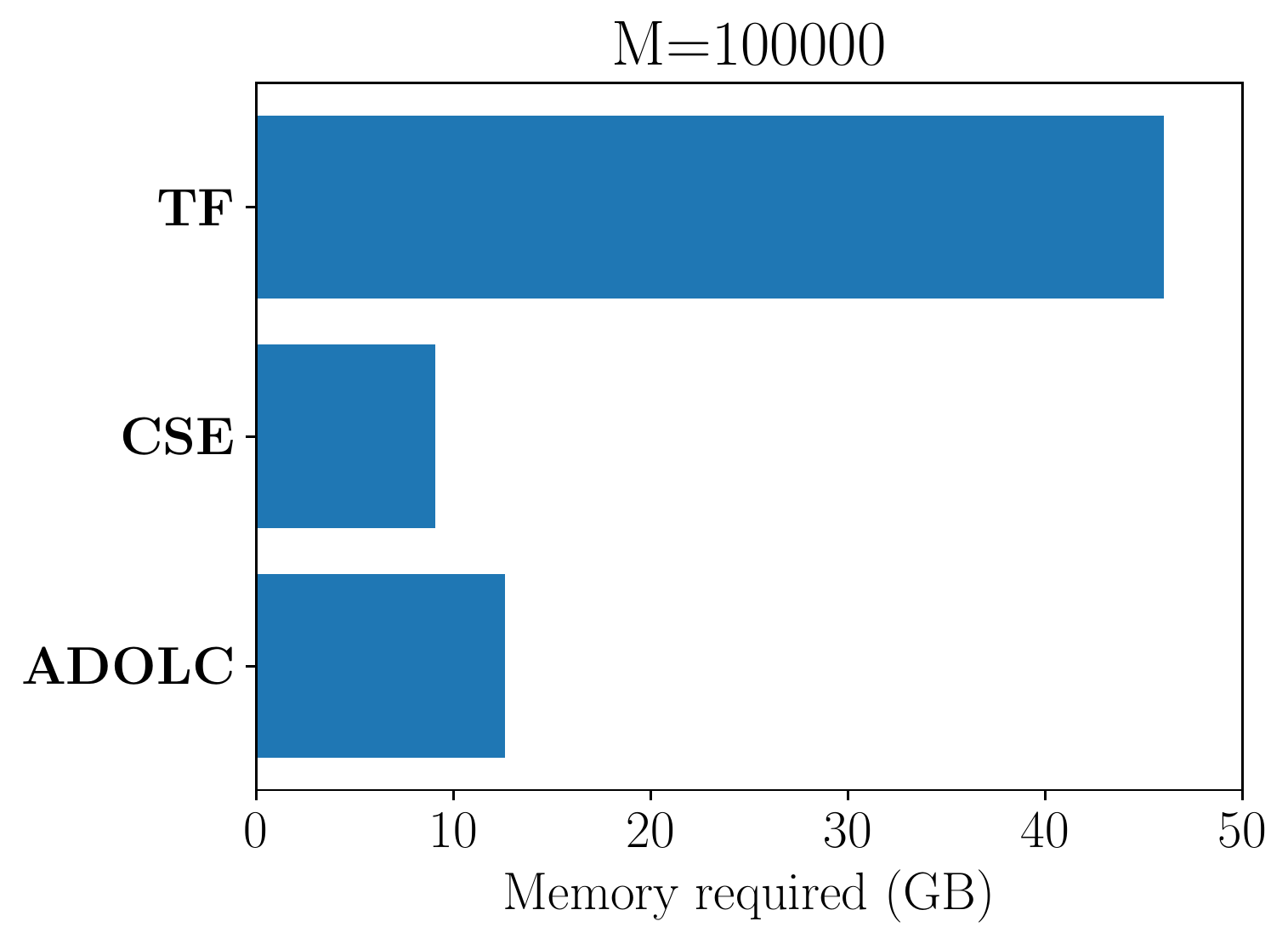}
	\end{subfigure}
	\caption{Memory consumption for different gradient calculation methods.  Note that Tensorflow (TF) has the highest memory footprint, 
which grows substantially with sample size. On the other hand, 
CSE, ADOLC and fnite-differences (FD) have fairly constant 
memory footprints. All numerical experiments were conducted on a 
dual Xeon E5-2651 v2 processor workstation operating 
at 1.80GHz with 145 GB of memory.}
	\label{fig:memory}
\end{figure}

It is difficult to quantify theoretically the memory and speed of 
graph-based methods to compute gradients such as 
\eqref{eqn:back_prop} . Hence, in this section we illustrate 
the memory savings and computational speedup we obtained 
through a series of numerical experiments. 
To this end, we consider the following model of wheel-driving 
UGV 
\begin{equation}
\dot{x}_1 = u_1 \cos(x_3),\qquad
\dot{x}_2 = u_1 \sin(x_3),\qquad 
\dot{x}_3 = \frac{u_1}{L} \tan(u_2),
\label{eq:UGV2} 
\end{equation}
where $(x_1(t),x_2(t))$ is the position of the vehicle, $x_3(t)$ 
is the heading angle, $u_1(t)$ is the forward velocity, and $u_2(t)$ is the steering angle (controls) -- see Figure \ref{fig:rear_car}(b). 
We set the following box constraints on the controls 
\begin{align*}
-1 \leq u_1(t) \leq 1,\qquad -1 \leq u_2(t) \leq 1\qquad \forall t\in[0,t_f].
\end{align*}
We integrate \eqref{eq:UGV2} in time from $t_0=0$ to $t_f=100$ using an explicit Runge Kutta 4th order method with step size 
$\Delta t=0.1$ (i.e., a total of $N=1000$ steps) 
and with a fixed, randomly generated, 
controls $u_1$ and $u_2$. Recall that in the multi-shooting 
setting we consider here, the controls are assumed 
to be constant over each given time interval 
$[t_i, t_{i+1}]$ (see Figure \ref{fig:directshooting}). 
To compute the gradient \eqref{eqn:back_prop}, hereafter we 
consider four different approaches, i.e., 
\begin{enumerate}
	\item Second-order centerd finite differences (FD)
	\item Operator overloading algorithmic differentiation (ADOLC \cite{walther2003adol})
	\item Graph-based algorithmic differentiation implemented in TensorFlow \cite{abadi2016tensorflow} (TF)
	\item Common Sub-expression Elimination (CSE)
\end{enumerate}
and compare them in terms of memory consumption, speed 
and accuracy. Owing to lack of an analytical solution, 
we will use ADOLC \cite{walther2003adol}, a mature algorithmic differentiation package to provide the benchmark gradient.
%
In this example we consider an objective that is the simple average 
of the UGV squared distance to the origin, namely
\begin{align}
J(\bm U) = \frac{1}{M} \sum_{i=1}^M \left[x_1^{(i)}(t_f)\right]^2
+\left[x_2^{(i)}(t_f)\right]^2,
\end{align}
where $M$ is the number of sample paths. 
In Figure \ref{fig:err1} we compare the absolute error 
in computing the gradient $\partial J/\partial u(t_j)$ for 
each $t_j$ in $[0,100]$.  It is seen that CSE and tensor flow 
are numerically exact, while second order (centered) 
finite-differences have accuracy of order $10^{-7}$,
corresponds to the threshold we set. 
Lower precision in gradient calculations typically results 
in more iterations during optimization when utilizing 
interior point methods. 
Next, we compare ADOLC, second-order finite-differences, tensorflow 
and CSE in terms of speed and memory consumption. This is done 
in Figure \ref{fig:runtime} and Figure \ref{fig:memory}, respectively. 
All numerical experiments were conducted 
on dual Xeon E5-2651 v2 processors operating at 1.80GHz with 145 GB of memory. 
For equivalence no GPUs were used.  ADOLC appears to 
scale best for very large samples,
but it is slower for samples below $M=10^5$
because of its lack of vectorization capabilities.
The numerical experiments we present\sout{ed} in this section
suggest that CSE and ADOLC and both suitable
algorithms for ensemble gradient calculations. CSE 
is faster and more memory efficient and can be 
vectorized on GPU devices, wheres ADOLC is currently 
restricted to run on CPUs.  
Tensorflow is extremely fast and suitable for solving only 
smaller problems because of its memory consumption 
and very long initialization times.

\section{Verification and validation of optimal controls}
\label{sec:V&V}
\begin{figure}[t!]
\centering
\includegraphics[width = 0.65 \textwidth]{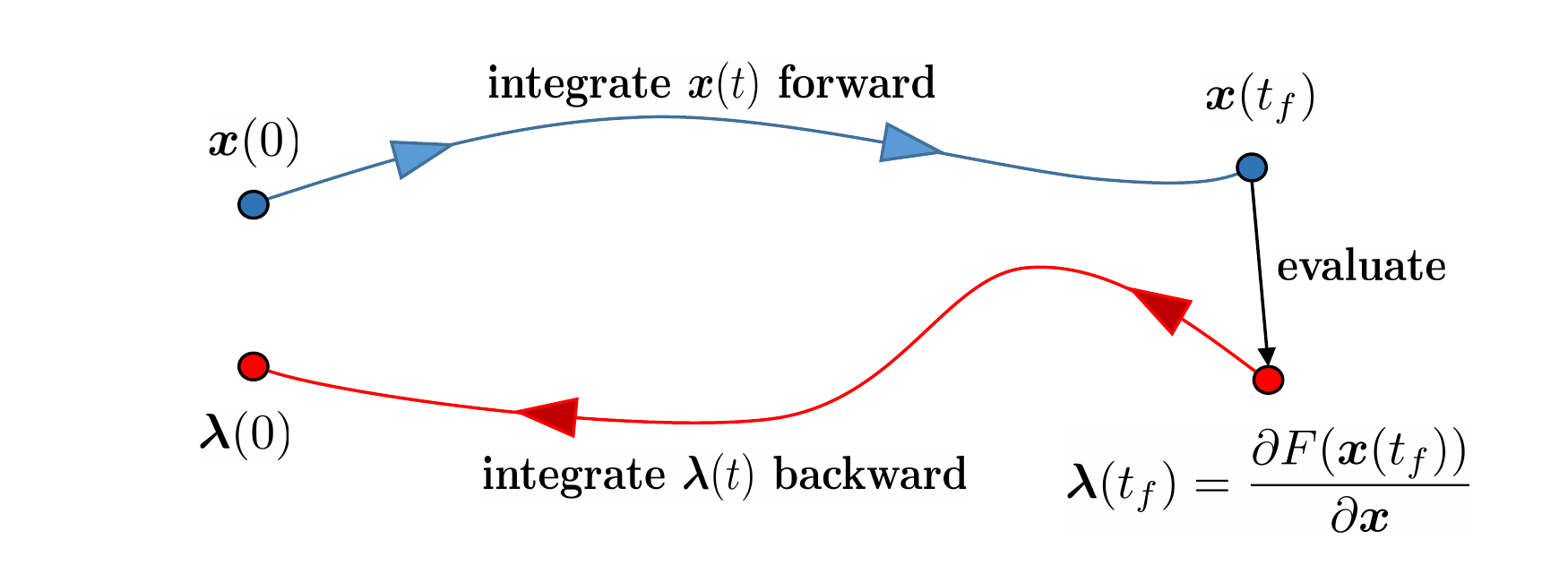}
\caption{Schematic of the forward-backward integration process for validation and verification of optimal controls.}
\label{fig: VandV diagram}
\end{figure}

In this section we propose a new criterion to 
verify and validate the numerical solution to the 
optimal control problem \eqref{eq:controlPr}. 
The criterion is based on the stochastic  Pontryagin's 
minimum principle \cite{JGCD_RS2015,Automatica_search}.
To describe method, consider the optimal 
control problem \eqref{eq:controlPr} 
and set  
\begin{equation}
J ([\bm x(t)],[\bm u(t)]) = \mathbb{E} \left \{ F (\bm x (t_f)) \right \} + 
\int_0^{t_f} r (\bm u(\tau)) d\tau,
\label{theJ}
\end{equation}
$r$, and $F$ being given functions. The second term at the right hand 
side is usually referred to as ``control energy'' term, and it can provide 
regularization of the optimization problem. 
Let us define the Hamilton's function
\begin{equation}
\label{eq: general Hamiltons function}
H (\bm u, \bm x, \bm \lambda) = r (\bm u) + 
\bm \lambda \cdot \bm f (\bm x, \bm u).
\end{equation}
From this we obtain the Hamilton's equations
\begin{equation}
\begin{dcases*}
\dot {\bm x} = \frac{\partial H}{\partial  \bm \lambda} = \bm f (\bm x, \bm u),\\
\dot {\bm \lambda} = - \frac{\partial H}{\partial \bm x} = - \frac{\partial \bm f}{\partial  \bm x} \cdot \bm \lambda.
\end{dcases*}
\label{eq: general Hamiltons equations}
\end{equation}
As shown in \cite{JGCD_RS2015}, 
the system \eqref{eq: general Hamiltons equations} 
is supplemented with a (random) initial condition for 
$\bm x(t)$, and a final condition for $\bm \lambda(t)$
\begin{equation}
\bm x (0) = \bm x_0 (\omega), \qquad 
\bm \lambda (t_f) = \frac{\partial F (\bm x (t_f))}{\partial \bm x}.
\label{infincond}
\end{equation}
In other words, \eqref{eq: general Hamiltons equations}-\eqref{infincond} 
is a two-point boundary value problem. 
Given any control $\widehat{\bm u}(t)$, we can solve \eqref{eq: general 
Hamiltons equations}-\eqref{infincond} using forward-backward 
propagation as illustrated in Figure \ref{fig: VandV diagram}. 
How can we check whether $\widehat{\bm u}(t)$ is optimal,
 i.e., it solves \eqref{eq:controlPr}-\eqref{theJ}?
A necessary condition is that the control satisfies the following 
minimization principle (extended Pontryagin's principle) 
\begin{equation}
\label{eq: general Hamilton true optimal}
\bm u^* (t) = \argmin_{\bm u (t)} \mathbb {E} \{H ({\bm u}(t), \bm x (t), \bm \lambda (t) \} \qquad \forall t \in [0, t_f],
\end{equation}
subject to the Hamilton's equations 
\eqref{eq: general Hamiltons equations}-\eqref{infincond}, and 
the constraints $\bm g(\bm u(t))\leq \bm 0$ on the control $\bm u(t)$. 
Note that both $\bm x(t)$ and $\bm \lambda(t)$ in 
\eqref{eq: general Hamilton true optimal} are functionals of $\bm u(t)$.
The Pontryagin's principle \eqref{eq: general Hamilton true optimal} 
offers us a simple way to verify the optimality of a given control.
To this end, suppose we have available a solution of 
\eqref{eq:controlPr}, say $\widehat{\bm u}(t)$,
computed using an optimizer, and we want to verify 
whether such solution is indeed optimal. A possible 
way to perform this calculation is to sample a large 
number of initial states $\bm x_0(\omega)$ and propagate 
them through the forward/backward problem 
\eqref{eq: general Hamiltons equations}-\eqref{infincond}. 
This allows us to obtain many realizations of $\bm x^{(i)} (t)$ and 
$\bm \lambda^{(i)} (t)$ for the given candidate 
control $\widehat{\bm u}(t)$. With such realizations available, 
we can approximate the expectation in 
\eqref{eq: general Hamilton true optimal} using, e.g., 
Monte Carlo and compute
\begin{equation}
\label{eq: general Hamilton false optimal}
\bm u^* (t) = \argmin_{\bm u (t)} \frac{1}{M}\sum_{i=1}^M 
H \left(\bm u, \bm x^{(i)} (t), \bm \lambda^{(i)} (t)\right).
\end{equation}
At this point we can compare $\bm u^* (t)$ 
with $\widehat{\bm u} (t)$ and verify whether 
the two controls are close to each other\footnote{We emphasize that the 
solution of \eqref{eq: general Hamilton false optimal} is {\em not}, 
in general, optimal since the expected value of the Hamiltonian 
is estimated using forward-backward propagation with 
fixed $\widehat{\bm u} (t)$. 
However, if $\left\| \widehat{\bm u} (t) - \bm u^* (t) \right\|$ 
is small then $\widehat {\bm u} (t)$ satisfies the necessary 
conditions for optimality.}.

\section{Numerical applications}
\label{sec:numerics}
In this section we demonstrate the proposed 
optimal control algorithms in applications to 
stochastic path planning problems involving models 
of an unmanned ground vehicles (UGVs) and
a fixed-wing unmanned aerial vehicles (UAVs).

\subsection{UGV stochastic path planning problem}
\label{sec:UGV}
\begin{figure}[t]
\centerline{
\includegraphics[height=10cm]{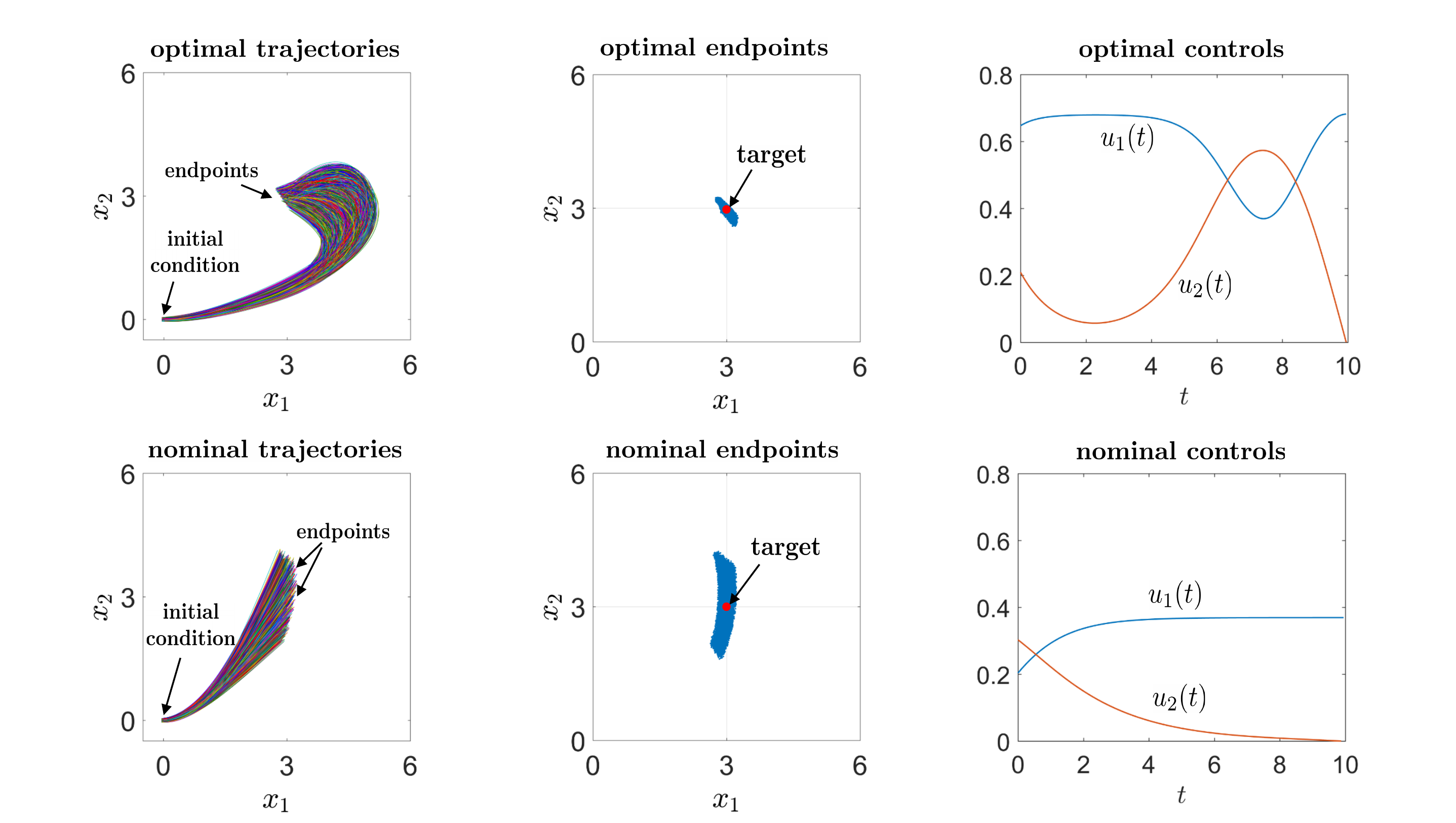}}
\caption{UGV stochastic path planning problem. Stochastic dynamics under 
optimal and nominal controls. It is seen that in both cases the 
endpoints of the UGV trajectories are clustered around 
the target located at $(x_1,x_2)=(3,3)$. However, in the 
optimal control case, to reduce the variance of the set of 
trajectories at final time, the UGV does not head directly towards 
the target, but it initially heads rights and then takes a steep left turn. 
On the other hand, in the nominal control case, the 
UGV heads directly towards the target, but the 
set of trajectories at final time is more spread out 
around the target than in the optimal control case.}
\label{fig:ugv_opt_traj}
\end{figure}

Consider the systems of equations \eqref{example:UGV}, describing 
the dynamics of a simple two-wheeled UGV with differential drive. We would like to control the velocity 
$u_1(t)$ and the steering angle $u_2(t)$ of the UGV in in a 
way that the vehicle hits a target located 
at $(x_1,x_2)=(3,3)$  at time $t_f=10$. To this end, we first 
assume that the initial state of the system is deterministic, i.e., 
\begin{equation}
x_{1}(0)=0, \qquad x_{2}(0)=0, \qquad x_{3}(0)=0,\qquad R=1.25. 
\label{initialstate_1}
\end{equation}
Moreover, we assume that the controls $u_1(t)$ and $u_2(t)$ 
are subject to the following constraints
\begin{equation}
-1 \leq u_1(t) \leq 1, \qquad 
-1 \leq u_2(t) \leq 1.
\label{UGV:control_constraints}
\end{equation}
We determined the optimal control for this deterministic 
problem by minimizing the objective functional
\begin{equation}
J([u_1(t),u_2(t)])= \frac{1}{2}\left\{\left(x_1(t_f)-3\right)^2+\left(x_2(t_f)-3\right)^2\right\}+\frac{q}{2}\int_{0}^{t_f}\left[u_1(t)^2+u_2(t)^2\right]dt.
\label{eqn:ugv-nom-obj}
\end{equation}
with respect to $u_1(t)$ and $u_2(t)$. Note that $x_1(t_f)$ and 
$x_2(t_f)$ are both functionals of $(u_1(t),u_2(t))$. 
We will refer to the controls minimizing \eqref{eqn:ugv-nom-obj} 
subject to the dynamics \eqref{example:UGV}, the initial condition 
\eqref{initialstate_1}  and the constraints 
\eqref{UGV:control_constraints} as the \textit{nominal controls}.  
This nominal control is used as the initial guess for the stochastic 
version of the problem, which can be formulated as follows: 
find controls $\{u_1(t),u_2(t)\}$ maximizing the probability 
that the UGV hits a target located at $(x_1,x_2)=(3,3)$ 
at $t_f=10$ under uncertain initial conditions and 
wheel radius $R$. In particular, we model such uncertainty 
as uniformly distributed independent random variables
\begin{equation}
x_{1}(0)\sim \mathcal U(-0.05, 0.05), \ \ 
x_{2}(0)\sim \mathcal U(-0.05, 0.05), \ \ 
x_{3}(0)\sim \mathcal U(-0.05, 0.05), \ \ 
 R\sim \mathcal U(1, 1.5).
 \label{eqn:ugv_unc}
\end{equation}
Note that the wheel radius $R$ here is assumed to be random. 
The stochastic version of the functional  \eqref{eqn:ugv-nom-obj}  
can be written as  
\begin{equation}
J([u_1(t),u_2(t)])=\frac{1}{2}\mathbb{E}
\left\{\left(x_1(t_f)-3\right)^2+\left(x_2(t_f)-3\right)^2\right\}+\frac{q}{2}\int_{0}^{t_f}\left[u_1(t)^2+u_2(t)^2\right]dt,
\label{JJ00}
\end{equation}
where $\mathbb{E}\{\cdot\}$ denotes an 
expectation over the (jointly uniform) 
PDF of $x_1(0)$, $x_2(0)$, $x_3(0)$, and $R$.
The optimal control problem reads as follows: 
minimize \eqref{JJ00} subject the dynamics
\eqref{UGV:control_constraints},  initial codnitions 
\eqref{eqn:ugv_unc} and control constraints \eqref{UGV:control_constraints}. 
The integral in both \eqref{eqn:ugv-nom-obj} and 
\eqref{JJ00} represents the energy of the control, 
and it provides a regularization of the control profiles, 
even for very small values of $q$. 
To compute the nominal and the optimal controls 
we utilized the multi-shooting algorithm we described 
in Section \ref{sec:multi-shooting}, with $S=2$ shooting 
intervals, explicit RK4 time integration ($\Delta t=0.05$), 
and CSE (section \ref{sec:CSE}) to calculate the ensemble gradient
of the cost functional (see Eq. \eqref{eqn:back_prop}). 
The number of sample paths to  approximate the expectation in \eqref{JJ00} with a Monte Carlo cubature were set to $M=10000$.
This yields a fully discrete optimal control problem of the form 
\eqref{eq:discrete_opt} (with no path constraints), 
which we solved using IPOPT \cite{Biegler,biegler2009large,xu2014pyipopt} 
linked with CSE and verified using Tensorflow (see Section \ref{sec:TensorFlow-IPOPT}). 
The results of our simulations are summarized in 
Figure \eqref{fig:ugv_opt_traj}. We see that, as expected, 
the final position of the UGV under nominal and optimal controls 
is clustered around the target located at $(x_1,x_2)=(3,3)$.  
However, in the optimal control case, to reduce the variance 
of the set of trajectories at final time, the UGV does 
not head directly towards the target, but it initially 
heads rights and then takes a steep left turn. 
On the other hand, in the nominal control case, the 
UGV heads directly towards the target, but the 
set of trajectories at final time is more spread out 
around the target than in the optimal control case.

\begin{figure}[t]
\centerline{
\includegraphics[height=6cm]{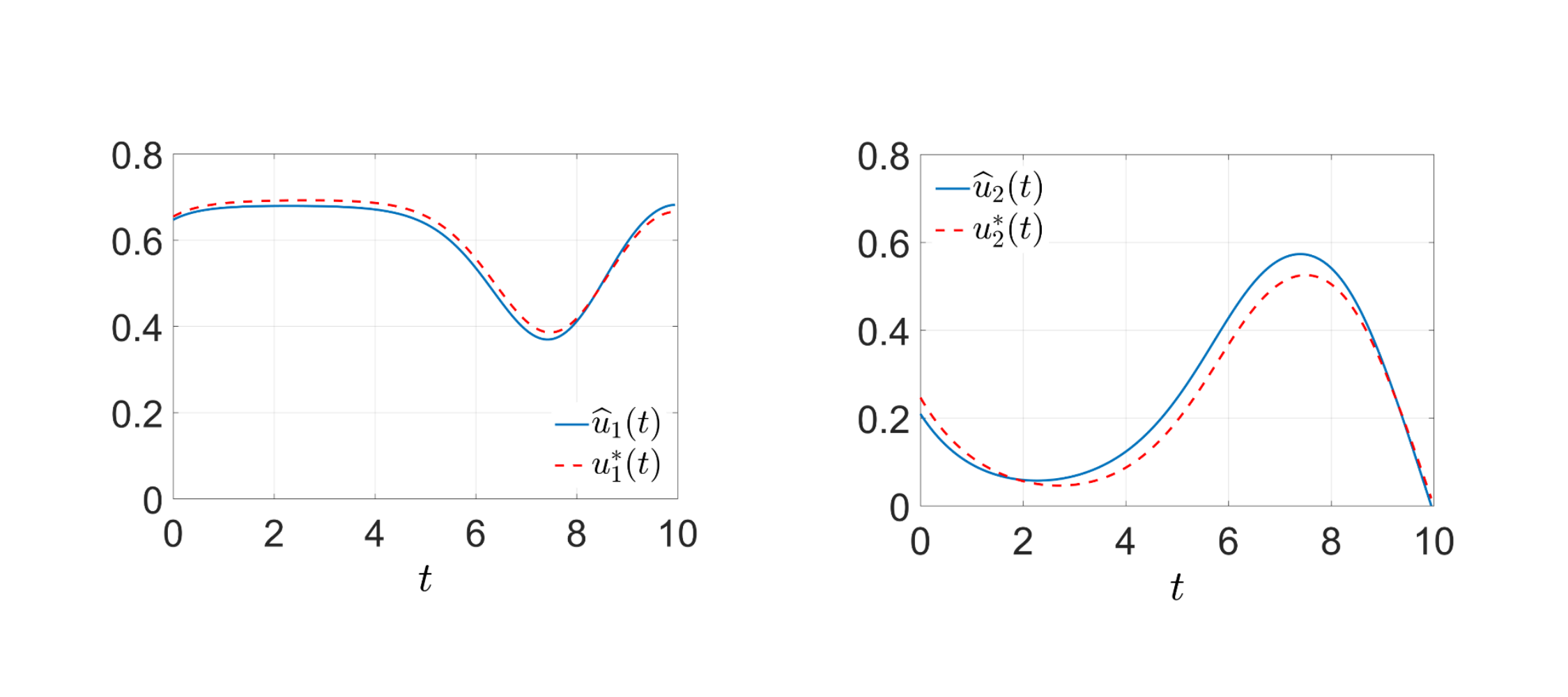}}
\caption{UGV stochastic path planning problem. Verification and 
Validation of the optimal control using the extended Pontryagin's 
minimum principle presented in Section \ref{sec:V&V}.}
\label{fig:V&V}
\end{figure}

For verification and validation of the optimal controls, 
we analytically derive the Hamiltonian function for the 
system \eqref{example:UGV}
\begin{equation}
\label{eq: segway Hamiltonian}
H (\bm u, \bm x, \bm \lambda) = \frac{q}{2} \left( u_1^2 + u_2^2 \right) + \lambda_1 R u_1 \cos (x_3) + \lambda_2 R u_1 \sin (x_3) + \lambda_3 R u_2 .
\end{equation}
The adjoint system is given by
\begin{equation}
\label{eq: segway adjoint}
\begin{dcases*}
\dot \lambda_1 = 0 , \\
\dot \lambda_2 = 0 , \\
\dot \lambda_3 = \lambda_1 R u_1(t) \sin (x_3) - \lambda_2 R u_1(t) \cos (x_3) .
\end{dcases*}
\end{equation}
We follow the process described in section \ref{sec:V&V} 
to measure the optimality of the computed solution.
This is done in Figure \ref{fig:V&V}, where we compare the control
$\bm u^*(t)=\{u_1^*(t),u_2^*(t)\}$ we obtained by solving 
\eqref{eq: general Hamilton false optimal} with 
optimal control $\widehat{\bm u}$ that minimizes 
\eqref{JJ00}. It is seen that the two controls are very close 
to each other, which means that $\widehat{\bm u}$ 
is a good approximation to the optimal control.

\subsection{UAV stochastic path planning problem}
\label{sec:UAV}
\begin{figure}[t]
\centering
\includegraphics[width=15cm]{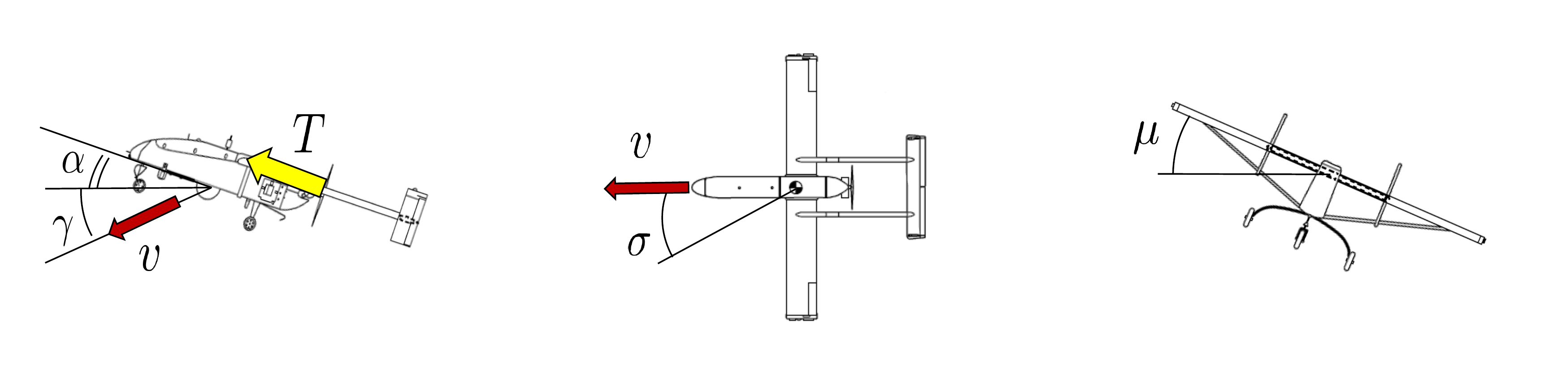}
\caption{Sketch of the fixed wing Unmanned Aerial Vehicle (UAV).
Shown are some of the phase variables appearing in the 
model \eqref{eq: UAV system}. The controls act on the 
trust $T$, the angle of attack $\alpha$, and the 
bank angle $\mu$.} \label{fig:UAVangles}
\end{figure}

In this section we consider a more challenging problem, namely 
stochastic path planning for a fixed-wing unmanned 
aerial vehicle (UAV) model under uncertain initial state 
such as position, heading angle, angle of attach and aerodynamic 
forces. The nonlinear dynamical system modeling the UAV 
is \cite{Savage2012,Shaffer2016}: 
	\begin{equation}
	\label{eq: UAV system}
	\begin{cases*}
	\dot x = v \cos \gamma \cos \sigma \\
	\dot y = v \cos \gamma \sin \sigma \\
	\dot z = v \sin \gamma \\
	\dot v = \displaystyle\frac{1}{m} (-D + T \cos \alpha) - g \sin \gamma \\
	\dot \gamma = \displaystyle\frac{1}{mv} (L \cos \mu + T \cos \mu \sin \alpha) - \frac{g}{v} \cos \gamma \\
	\dot \sigma = \displaystyle\frac{1}{mv \cos \gamma} (L \sin \mu + T \sin \mu \sin \alpha) \\
	\dot T = u_T \\
	\dot \alpha = u_{\alpha} \\
	\dot \mu = u_{\mu}
	\end{cases*}
	\end{equation}
where $(x(t),y(t),z(t))$ is the position of the UAV in a Cartesian reference 
frame, $v(t)$ is the velocity, $(\gamma(t), \sigma(t))$ 
are elevation and heading angles,  $T(t)$ is the thrust,  $\alpha(t)$ 
is the angle of attack, and $\mu(t)$ is the bank angle 
(see Figure \ref{fig:UAVangles}).  The UAV model \eqref{eq: UAV system} 
is valid only in the region of the phase space defined by the following 
constraints
\begin{align}
13 \le v(t) \le 42, \quad
-\frac{\pi}{6} \le \gamma(t) \le \frac{\pi}{6}, \quad 
-\pi \le \sigma(t) \le \pi, \quad 
3.0 \le T(t) \le 35.0, \quad 
-\frac{\pi}{8} \le \alpha(t) \le \frac{\pi}{8}.
\label{UAVstatespace}
\end{align}
These are effectively linear state space constraints 
that needs to be added to the optimal control problem 
(see Eqs. \eqref{eq:controlPr} and \eqref{eq:discrete_opt}).
The other parameters appearing in \eqref{eq: UAV system} 
are the UAV mass $m=2$ kg, the acceleration of gravity 
$g = 9.8$ m/s$^2$, and the aerodynamic lift and drag 
forces given by
\begin{align}
\label{eq: lift and drag}
L = \frac{1}{2} \rho v^2 S C_L,
\qquad
D = \frac{1}{2} \rho v^2 S C_D, 
\end{align}
In these expressions $\rho(z) = 1.21 e^{-z / 8000}$ kg/m$^3$ 
is the mass density of air, 
$S = 0.982$ m$^2$ is the wing surface area, and $C_L$ and $C_D$ 
are the lift and drag coefficients defined as
\begin{align}
C_L &= (C_{x0} + C_{xa} \alpha) \sin \alpha - 
(C_{z0} + C_{za} \alpha) \cos \alpha,
\label{eq: lift coefficient}\\
C_D &= -(C_{x0} + C_{xa} \alpha) \cos \alpha - 
(C_{z0} + C_{za} \alpha) \sin \alpha.
\label{eq: drag coefficient}
\end{align}
Note that $C_L$ and $C_D$ depend on the parameters $C_{x0}, C_{xa}, C_{z0},$ and $C_{za}$ which are usually unknown and 
must be estimated from data. In our simulation we model 
such coefficients as independent random variables with 
given probability distributions. 

\begin{figure}[t]
		\centering
		\begin{subfigure}{.5\textwidth}
			\centering
			\includegraphics[width=.9\linewidth]{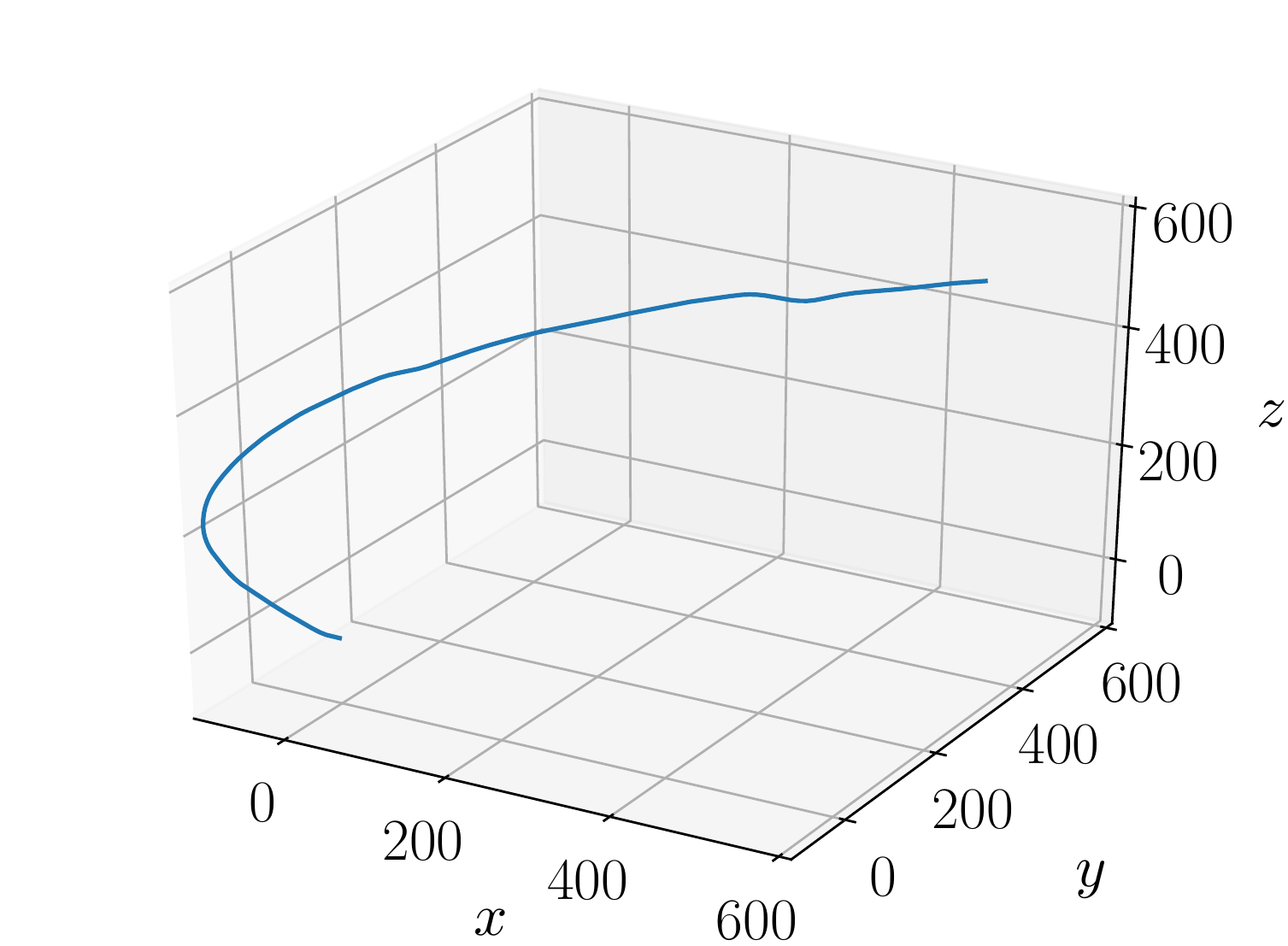}
			\caption{\footnotesize Nominal trajectory.}
			\label{fig:uav_nom_traj.pdf}
		\end{subfigure}%
		
		\begin{subfigure}{.5\textwidth}
			\centering
			\includegraphics[width=.9\linewidth]{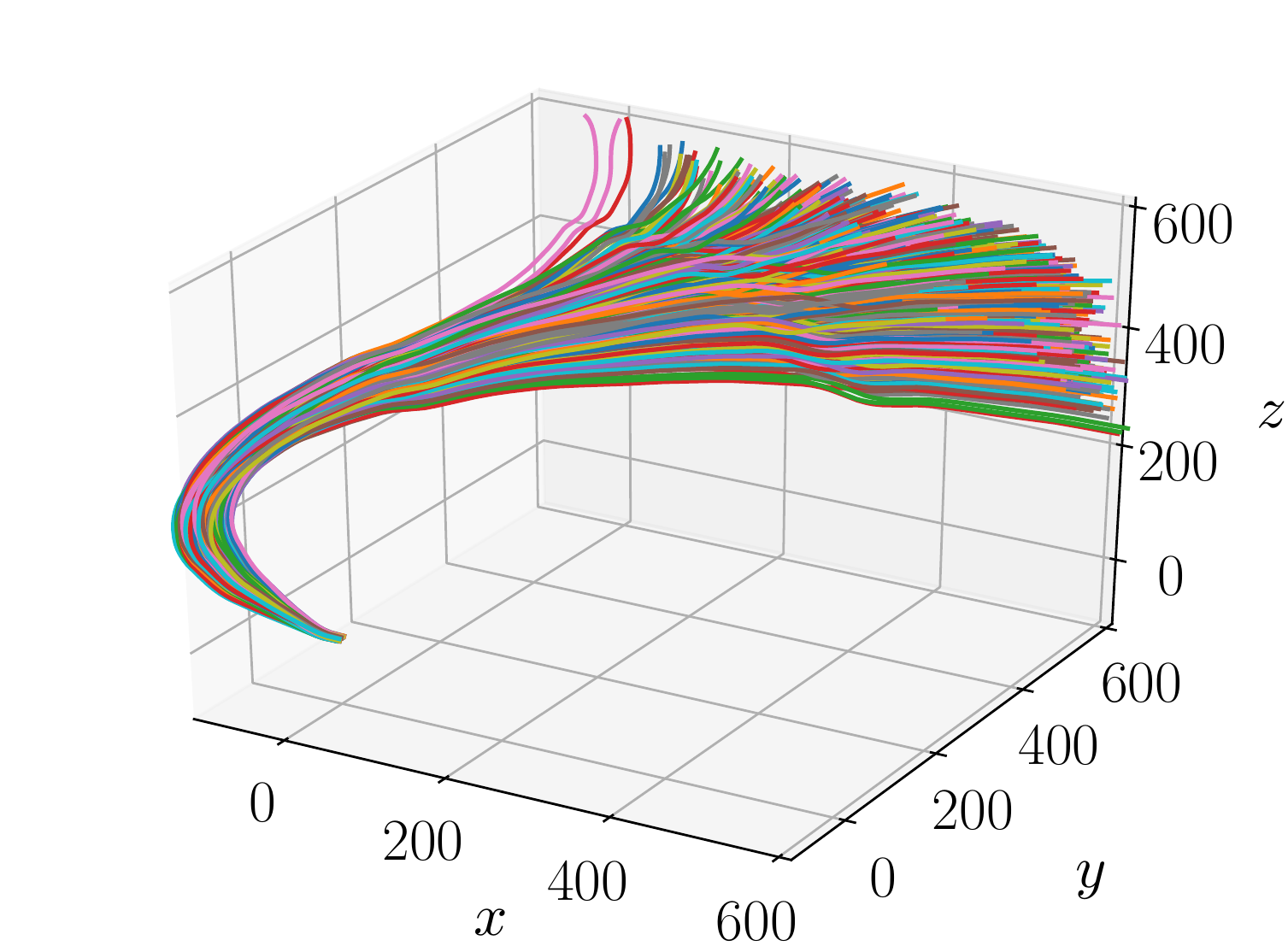}
			\caption{\footnotesize  Random trajectories under nominal controls.}
			\label{fig:uav_nom_unc}
		\end{subfigure}%
		\begin{subfigure}{.5\textwidth}
			\centering
			\includegraphics[width=.9\linewidth]{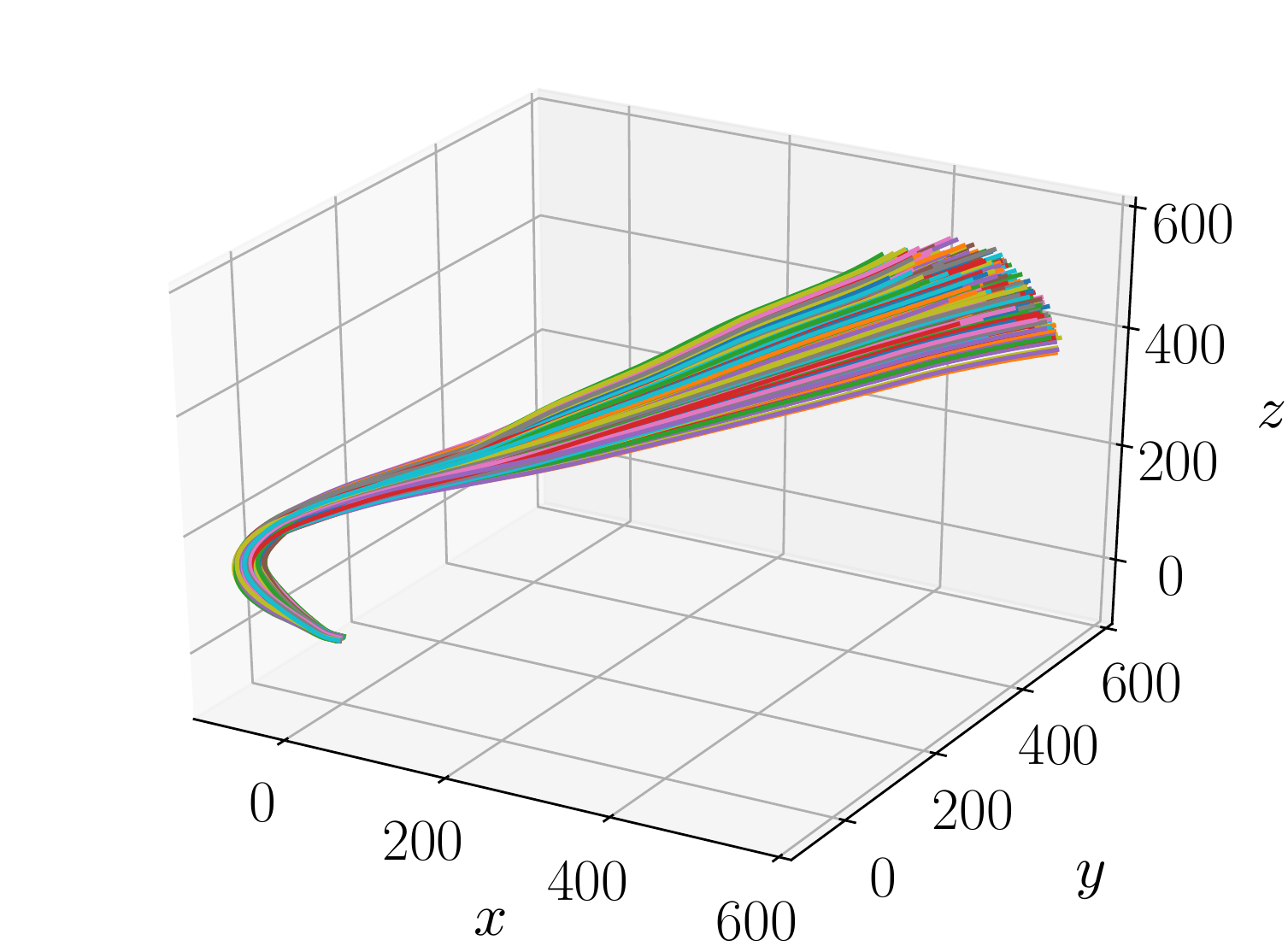}
			\caption{\footnotesize  Random trajectories under optimal controls.}
			\label{fig:uav_nom_T}
\end{subfigure}%
\caption{\footnotesize UAV stochastic path planning problem. 
Shown are trajectories obtained by sampling the random 
initial condition and the random parameters of the model 
\eqref{eq: UAV system} and then propagating them forward in time 
using nominal and optimal controls. It is seen that the optimal controls 
can mitigate the effects of uncertainty, and they yield 
final states clustered around the target located at 
$(x,y,z)=(500,500,500)$.  Figures drawn to the same scale 
for comparison.}
\label{fig:uav_trajs}
\end{figure}

We aim at computing optimal controls for the angle of attach $u_\alpha(t)$, 
the bank angle $u_\mu(t)$ and trust $u_T(t)$ so that the UAV hits 
the target located at $(x,y,z)=(500,500,500)$ 
under uncertain aerodynamic forces, uncertain initial position/angles, 
and uncertain initial velocity.
The controls are subject to the following box constraints.
\begin{align}
u_T \in [-1.0, 1.0],\qquad 
u_\alpha \in [-0.05, 0.05],\qquad
u_\mu \in [-0.05, 0.05].
\label{UAVconst}
\end{align}
Following the same steps as in the UGV example we studied 
in section \ref{sec:UGV}, we first calculate the nominal control 
where we minimize the functional 
\begin{align}
J\left( [u_T(t)], [u_{\alpha}(t)], u_{\mu}(t)]\right)=&
(x(t_f)-500)^2+(y(t_f)-500)^2+(z(t_f)-500)^2 + \nonumber\\
&\frac{1}{200}\int_0^{t_f}
\left[u_T^2(t)+u_{\alpha}^2(t)+u_{\mu}^2(t)\right]dt
\end{align}
subject to \eqref{eq: UAV system}, \eqref{UAVstatespace}, 
\eqref{UAVconst} and the deterministic initial states
\begin{gather*}
x(0) = 0, \quad 
y(0) = 0, \quad 
z(0) = 0, \quad 
v(0) = 27.5,\quad
\gamma(0) = 0,\quad 
\sigma(0) = \pi,\quad
T(0) = 16.1,\nonumber\\
\alpha(0) =  -0.0088,\quad 
\mu(0) =0,\quad 
C_{x0} = -0.03554,\quad 
C_{xa} = 0.00292,\quad  
C_{z0}= -0.055,\quad 
C_{za}= -5.578.\nonumber
\end{gather*}
Note that this defines a challenging 
maneuver whereby the UAV has to move from  
$(x(0),y(0),z(0))=(0,0,0)$ to $(x(t_f),y(t_f),z(t_f))=
(500, 500, 500)$ while the initial heading 
angle $\sigma$ is in the opposite direction 
(see Figure \ref{fig:uav_nom_traj.pdf}).

Next, we consider the stochastic path planning problem, i.e., 
the problem of computing robust controls that 
can steer the UAV from an uncertain initial state 
(position, velocity, and configuration angles) 
to the final position $(x(t_f),y(t_f),z(t_f))=
(500, 500, 500)$, under random aerodynamics forces.
A simplified version of this problem was recently studied 
by Shaffer {\em et. al} in \cite {Shaffer2016}, where  
uncertainty was modeled in terms of only 
one random variable, i.e., $C_{x0}$ in Eqs. 
\eqref{eq: lift coefficient} and  \eqref{eq: drag coefficient}. 
The main reason for studying such simplified model was 
that the computational control algorithm could not 
handle higher-dimensional random input vectors. 
Indeed, the memory requirements and computational cost 
of the algorithms proposed in \cite{Shaffer2016}
are currently beyond the capabilities of a modern 
workstation.  
To compute the optimal controls in the stochastic path planning 
UAV problem we minimize the cost functional 
\begin{align}
\label{eqn:UAV_nom}
J\left([u_T(t)], [u_{\alpha}(t)], [u_{\mu}(t)]\right) = &\mathbb{E}\left\{(x(t_f)-500)^2+(y(t_f)-500)^2+(z(t_f)-500)^2\}\right\} + \nonumber\\
&\frac{1}{200}\int_0^{t_f}\left[u_T^2(t)+u_{\alpha}^2(t)+u_{\mu}^2(t)\right]dt,
\end{align}
subject to \eqref{eq: UAV system}, \eqref{UAVstatespace}, 
the random initial condition 
\begin{gather}
x(0) = r\sin(\theta) \cos(\phi), \qquad 
y(0) = r\sin(\theta) \sin(\phi), \qquad 
z(0) = r\cos(\phi),  \label{eq: UAV initial PDFs0}\\
v(0) \sim  \mathcal U (25.5750, 29.4250), \qquad
\gamma(0)\sim  \mathcal{U} (-0.05, 0.05),\qquad
\sigma(0) \sim  \mathcal{U} (3.1, 3.2),\label{eq: UAV initial PDFs1}\\ 
C_{x0} \sim  \mathcal U  (-0.0380, -0.0330),\qquad
C_{xa} \sim  \mathcal U (0.0027, 0.0031),\nonumber\\
C_{z0} \sim  \mathcal U (-0.0589, -0.0548),\qquad  
C_{za} \sim  \mathcal U (-5.9685, -5.1875),
\label{eq: UAV initial PDFs2}
\end{gather}
where $r \sim \mathcal U(0, 5.0)$, $\theta \sim \mathcal U(0, \pi)$ 
$\phi \sim \mathcal U(0, 2 \pi)$, and the ensemble path constraints 
\begin{gather*}
13 \le \mathop{\mathbb{E}} \{ v(t) \} \le 42, \qquad
-\frac{\pi}{6} \le \mathop{\mathbb{E}} \{ \gamma(t) \} \le \frac{\pi}{6},\qquad 
-\pi \le \mathop{\mathbb{E}} \{ \sigma(t) \} \le \pi, \\ 
3.0 \le \mathop{\mathbb{E}} \{T(t)\} \le 35.0,\qquad 
-\frac{\pi}{12} \le \mathop{\mathbb{E}} \{\alpha(t)\} \le \frac{\pi}{12}.
\end{gather*}
In Eqs. \eqref{eqn:UAV_nom} and \eqref{eq: UAV system} 
$\mathbb{E}\{\cdot \}$ is an expectation over the joint PDF of 
the random variables in the initial state, which we approximate using 
Monte Carlo cubature rule with $M=4800$ randomly drawn sample 
paths. 
\begin{figure}[t!]
\centerline{
\includegraphics[width=5.5cm]{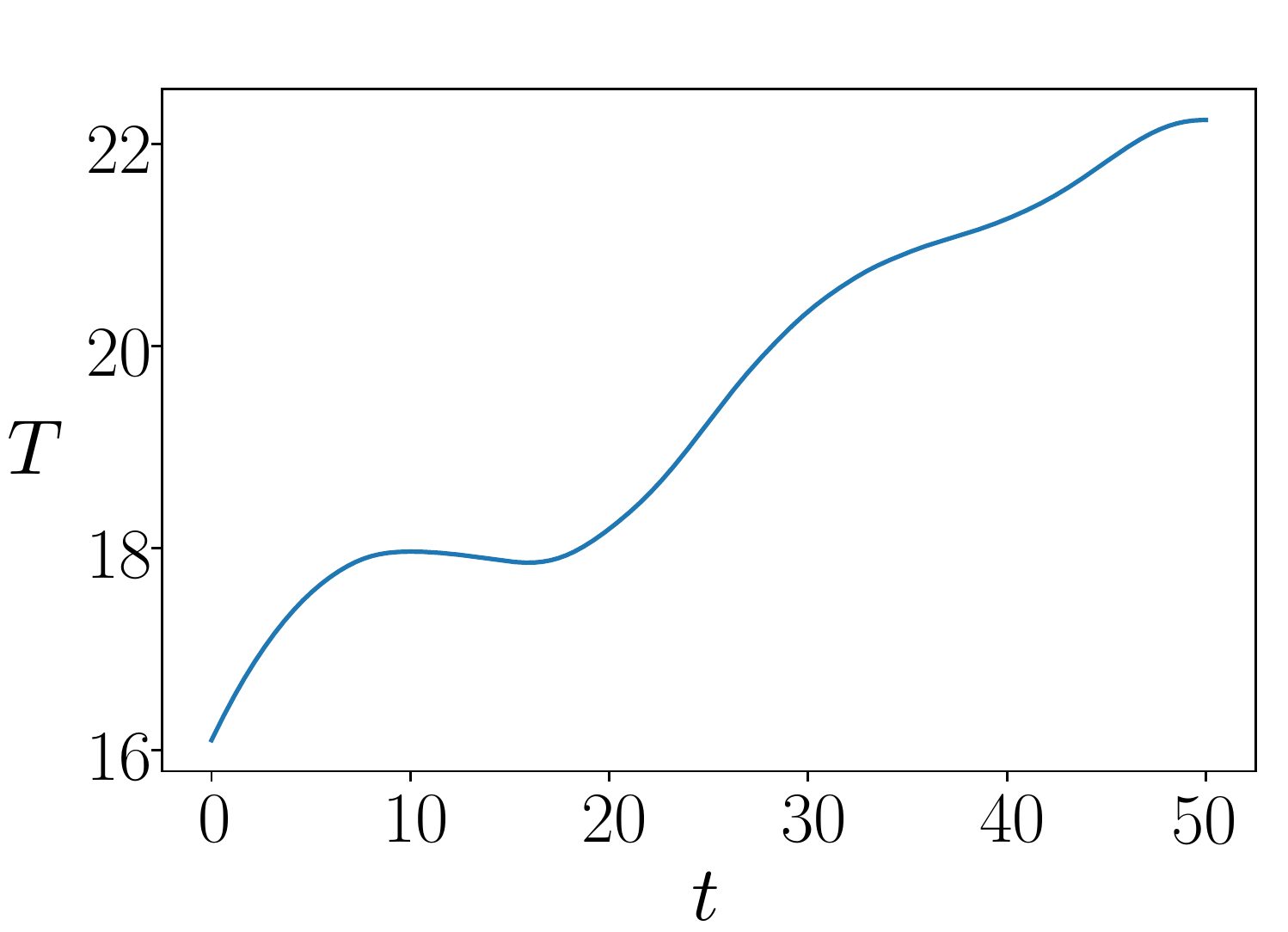}
\includegraphics[width=5.2cm]{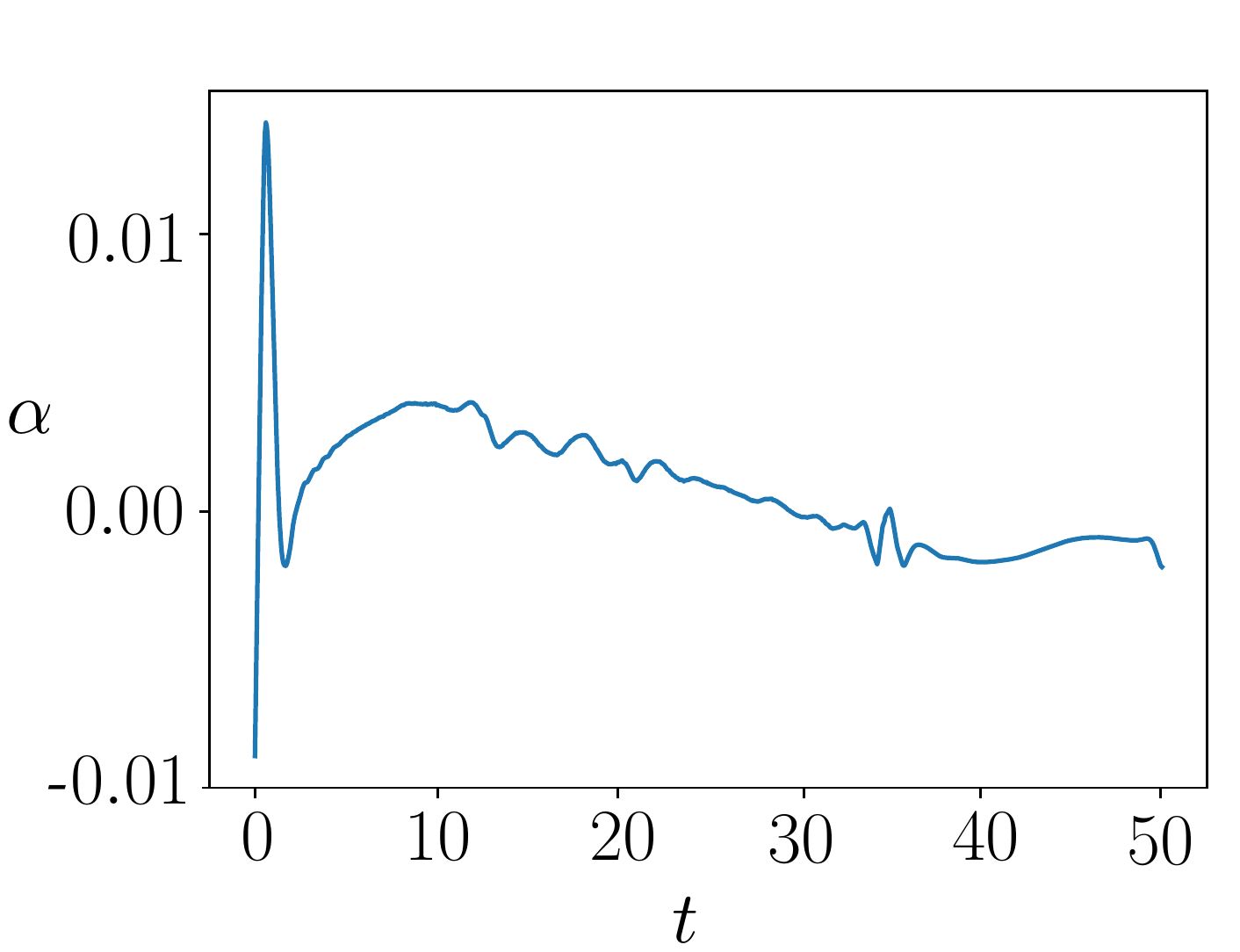}
\includegraphics[width=5.5cm]{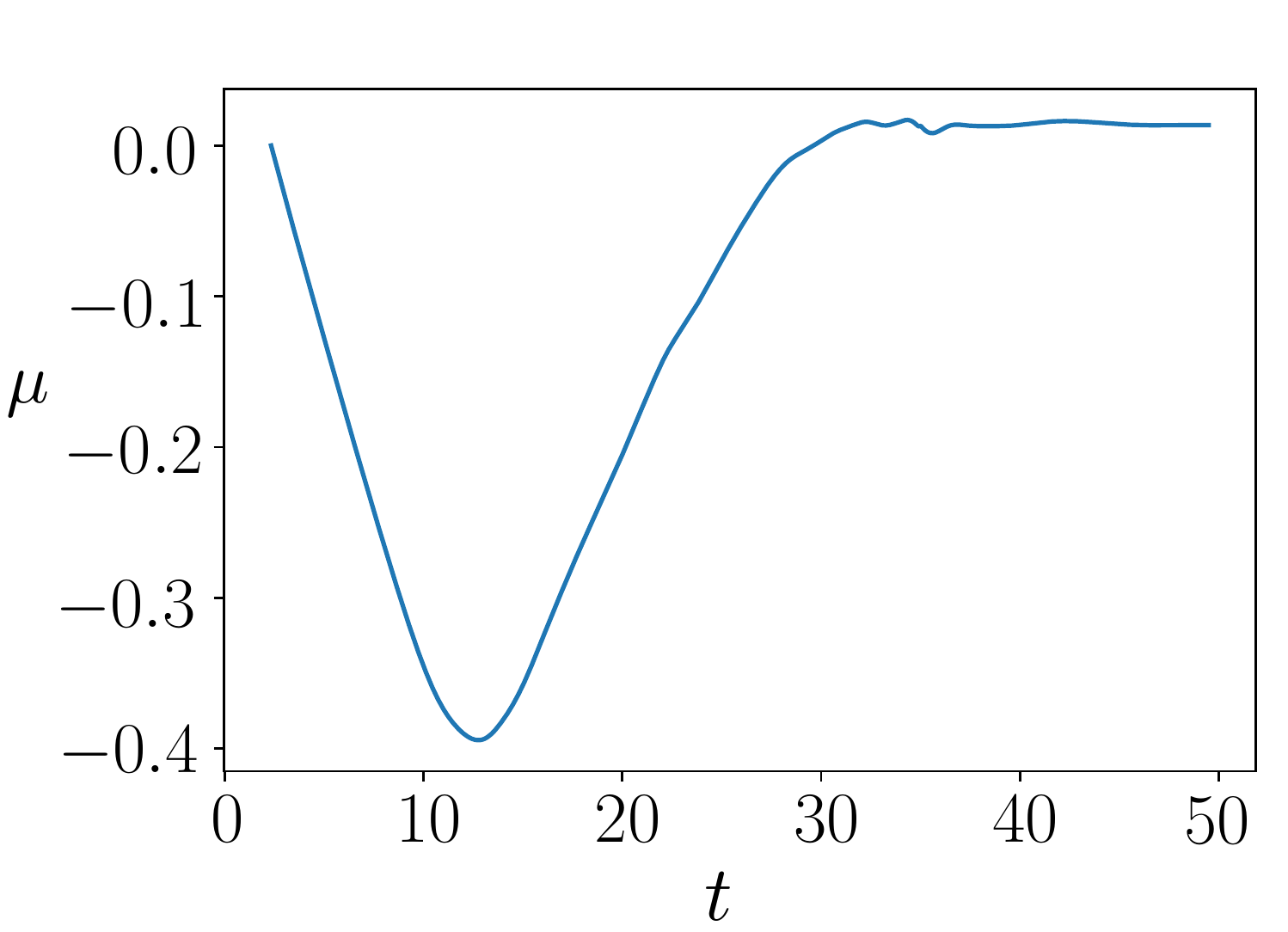}}
\caption{Stochastic path planning problem for the 
fixed-wing UAV model \eqref{eq: UAV system}. 
Shown are the optimal trust $T(t)$, angle of attack 
$\alpha(t)$ and bank angle $\mu(t)$ 
(integrals of $u_T$, $u_\alpha$ and $u_\mu$) 
we obtained to steer the UAV from the random initial 
position \eqref{eq: UAV initial PDFs0} to the final position 
$(x,y,z)=(500,500,500)$, under random 
aerodynamic forces.}
\label{fig:uav_opt_u}
\end{figure}
To compute the nominal and the optimal controls 
we utilized a single shooting algorithm 
with explicit RK3 time integration ($\Delta t=0.1$), 
and CSE to calculate the ensemble gradient 
of the cost functional (see Eq. \eqref{eqn:back_prop}). 
This yields a fully discrete optimal control problem 
of the form \eqref{eq:discrete_opt}, 
which we solved using IPOPT \cite{Biegler,biegler2009large,xu2014pyipopt} 
linked with CSE, and verified using ADOLC  \cite{walther2003adol}. 
The optimal  trust, angle of attack and bank angle 
are shown in Figure \ref{fig:uav_opt_u}. To verify the 
performance of the optimal controls under random 
initial conditions and uncertain parameters, we propagated in time 
the controlled dynamics of $10000$ randomly generated
trajectories, with initial conditions and parameters sampled
according to  \eqref{eq: UAV initial PDFs0}-\eqref{eq: UAV initial PDFs2}. 
Such trajectories are shown in Figure \ref{fig:uav_trajs}.  
It is seen that the final positions of the UAV converge 
to a smaller region about the target 
location $(x,y,z)=(500,500,500)$. 

\subsubsection{Nonlinear advection-reaction-diffusion PDE}

\begin{figure}[t]
\centerline{\includegraphics[height = 5.5cm]{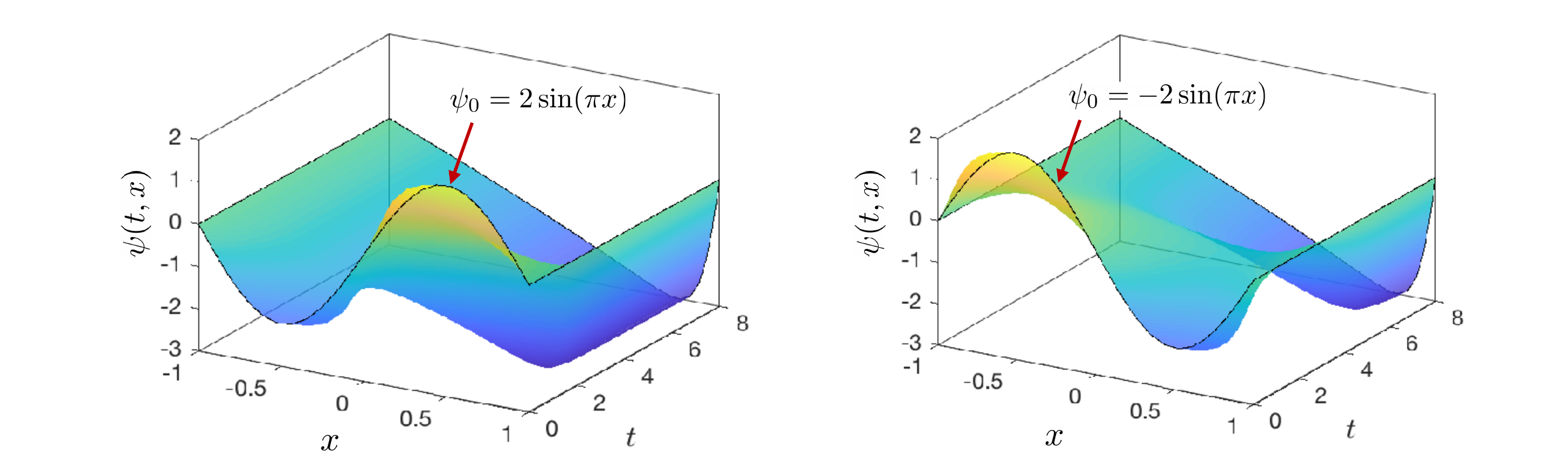}}
\caption{Solution samples of the initial-boundary value problem 
\eqref{eq: Burgers PDE} corresponding to different initial conditions. 
Here we set the control $u(t)$ equal to zero (uncontrolled dynamics).}
\label{fig:PDEsamples}
\end{figure}

In this section we study open-loop control under uncertainty 
of nonlinear PDEs. Specifically, we consider the following 
initial-boundary value problem involving a nonlinear 
advection-reaction-diffusion equation 
\begin{equation}
\label{eq: Burgers PDE}
\left \{ \begin{array}{ll}
\displaystyle \frac{\partial \psi}{\partial t} = \psi\frac{\partial \psi}{\partial x} + \frac{1}{5}
\frac{\partial^2 \psi}{\partial x^2} +\frac{3}{2} \psi e^{-\psi/10} + 
I_{\Omega} (x) u(t), \vs\\
\psi (t, -1) = \psi(t, 1) = 0,  \vs\\
\psi (0, x) = \psi_0(x),
	\end{array} \right.
\end{equation}
where $\psi(x,t): [-1,1]\times [0,t_f]\mapsto \mathbb{R}$ is a random 
field that is a functional of $\psi_0$ (random initial condition) and  
$u (t): [0, t_f] \to \mathbb R$ (control). Note that $u(t)$  is actuated 
only on the spatial domain $\Omega=[-0.5,-0.2]$, which is 
the support of the indicator function $I_{\Omega} (x)$.
We aim at determining the optimal control $u(t)$ that minimizes 
the cost functional 
\begin{equation}
J([\psi],[u])=\frac{1}{2}\mathbb{E}\left\{\int_{-1}^1\psi(t_f,x)^2dx +2\int_{0}^{t_f}
\int_{-1}^1\psi(\tau,x)^2dxd\tau \right\}+\frac{1}{20}\int_{0}^{t_f}u(\tau)^2d\tau,
\qquad t_f=8,
\label{eq:costPDE}
\end{equation}
where $\mathbb{E}\{\cdot\}$ is an expectation over 
the probability distribution of the initial state $\psi_0(x)$. 
To this end, we first discretize \eqref{eq: Burgers PDE} in space using 
a Chebyshev pseudo-spectral method \cite{Trefethen2000}. 
This yields the semi-discrete form
\begin{equation}
\left \{ \begin{array}{ll}
\displaystyle \dot{\bm \psi} = \bm \psi \circ \bm {D}\bm \psi + \frac{1}{5} 
{\bm D}^2 \bm \psi + \frac{3}{2} \bm \psi\circ  e^{-\bm \psi/10}  + 
\bm I_{\Omega} u(t), \vs\\
\bm \psi(0)= \bm \psi_0,
\end{array}\right.
\label{eq:semidiscrete}
\end{equation}
where $\bm \psi(t)=[\psi(t,x_1),...,\psi(t,x_{n})]^T$ collects 
the values of the solution $\psi(t,x)$ at the (inner) Chebyshev 
nodes
\begin{equation}
x_j = \cos \left(\frac{j \pi}{n+1}\right)\qquad j=1,...,n.
\end{equation}
In equation \eqref{eq:semidiscrete}, the circle ``$\circ$'' 
denotes element-wise multiplication (Hadamard product), 
$\bm I_{\Omega}$ is the discrete indicator function, 
while $\bm D$ and $\bm D^2 \in \mathbb R^{n \times n}$ are, 
respectively, the first- and second-order Chebyshev differentiation 
matrices. These matrices obtained by deleting the first 
and last rows and columns of the full differentiation matrices. 
In Figure \ref{fig:PDEsamples} we plot two solution samples 
of the initial-boundary value problem \eqref{eq: Burgers PDE} 
we obtained by integrating the discretized 
system \eqref{eq:semidiscrete} in time.
The cost functional \eqref{eq:costPDE} can be discretized as 
\begin{equation}
J = \frac{1}{2M}\sum_{i=1}^M \bm w^T \left[ \bm \psi^{(i)}(t_f)^2+
2\int_{0}^{t_f}\bm \psi^{(i)}(\tau)^2d\tau\right]+\frac{1}{20}
\int_{0}^{t_f}u(\tau)^2d\tau,
\label{discretefunct}
\end{equation}
where $\bm w$ are Clenshaw-Curtis quadrature 
weights (column vector), and $M$ is the 
total number of samples. 

\begin{figure}[t]
\centerline{\includegraphics[height=5.5cm]{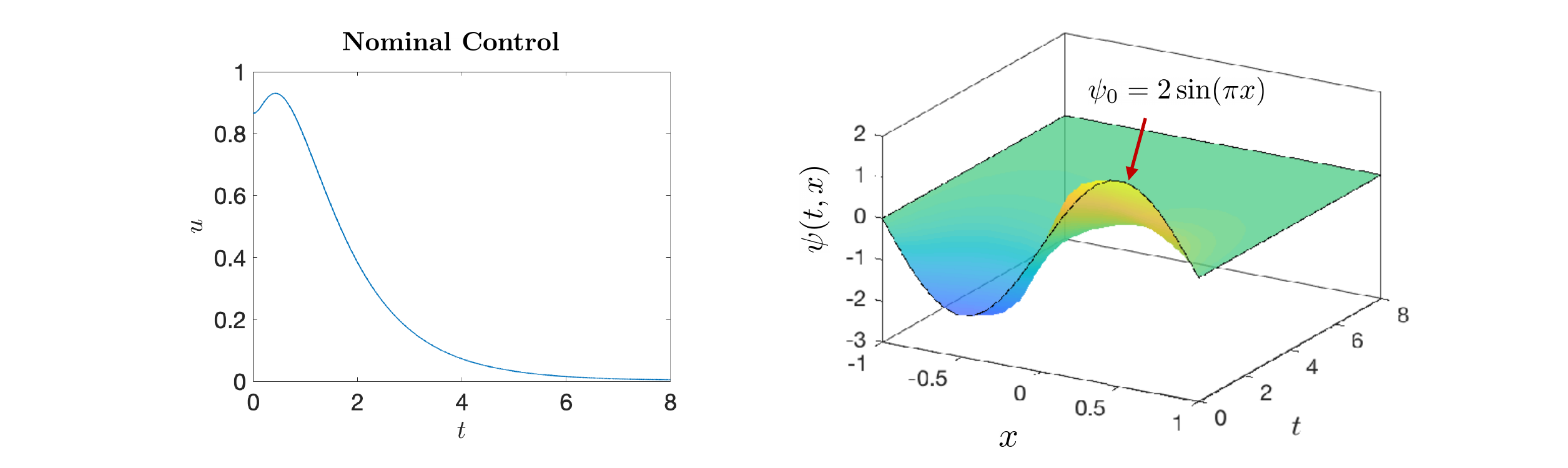}}
\caption{Open loop control of the nonlinear PDE \eqref{eq: Burgers PDE}. 
Shown is the nominal control minimizing \eqref{eq:costPDE} for 
the deterministic initial condition $\psi_0=2\sin(\pi x)$ and the 
corresponding solution dynamics. It is seen 
that the control $u(t)$ sends $\psi(t,x)$ to zero 
after a small transient.}
\label{fig:controlledBurgers}
\end{figure}

\begin{figure}[t]
\centerline{\includegraphics[height=9.8cm]{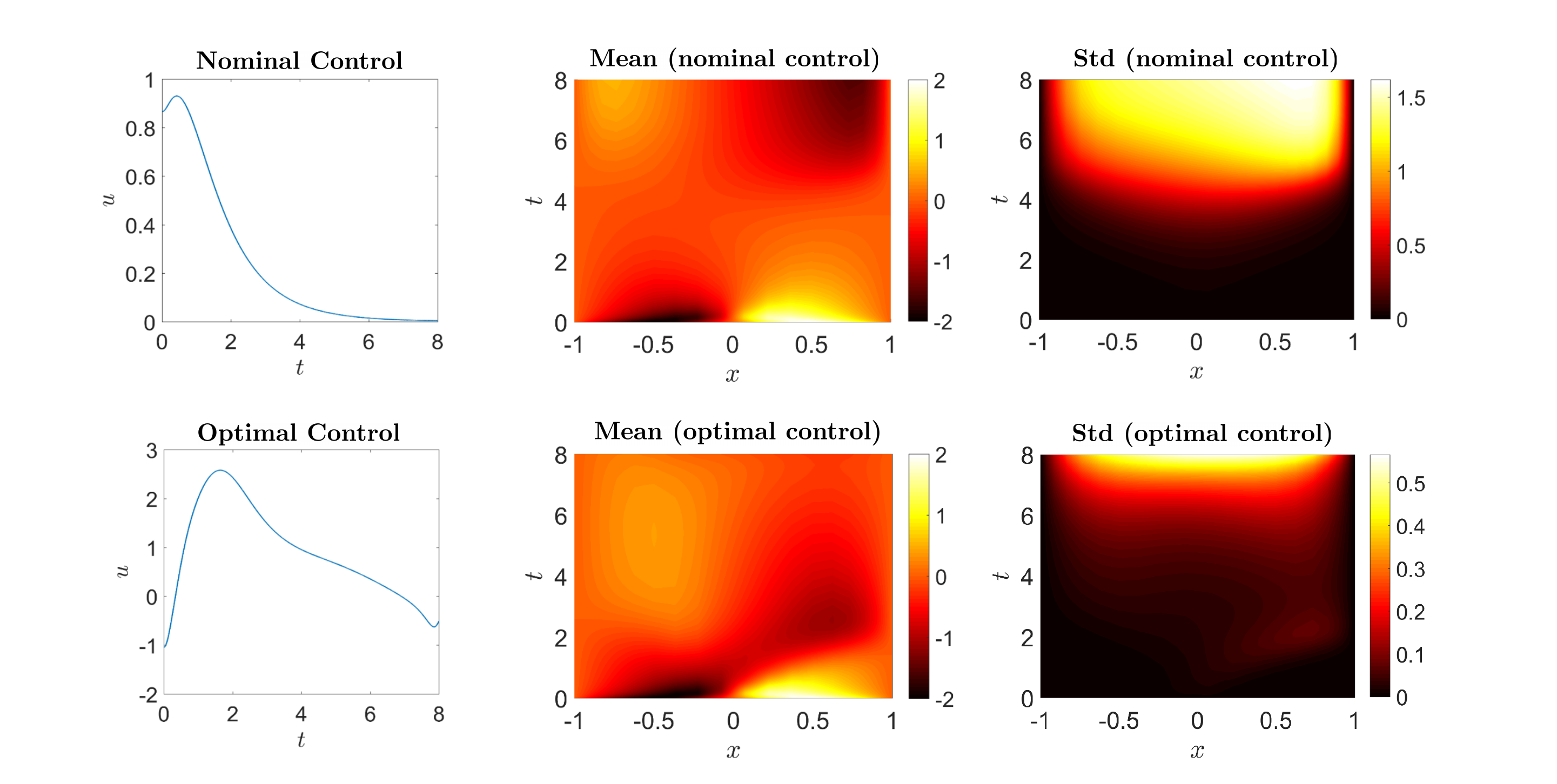}}
\caption{Comparison between the mean and the standard deviation of the 
stochastic solution to the PDE \eqref{eq: Burgers PDE} 
under optimal or nominal controls. The initial condition here is set to be 
random as \eqref{eq:randomIC}. It is seen that the optimal control 
is more effective is driving the stochastic solution to zero. 
In fact, both the mean and standard deviation of the 
optimally controlled solution are closer to zero. 
Note also that the optimal control differs 
substantially from the nominal control.}
\label{fig:PDEcontrolcomparison}
\end{figure}

To quantify the effectiveness of the ensemble optimal control algorithm 
we propose, we first determine the nominal control corresponding to the 
deterministic initial condition $\psi_0(x)=2\sin(\pi x)$. Such control minimizes
the functional \eqref{discretefunct} with $M=1$ and $\psi^{(1)}_j(0)=2\sin(\pi x_j)$, 
and sends $\psi(t,x)$ to zero after a small transient (see Figure 
\ref{fig:controlledBurgers}). 
Next, we introduce uncertainty in the initial condition. 
Specifically, we set 
\begin{equation}
\bm \psi_{0j} = 2\sin(\pi x_j)+\epsilon_j\qquad  \epsilon_j\sim \mathcal{N}(0, 001). 
\label{eq:randomIC}
\end{equation}
This makes the solution to \eqref{eq: Burgers PDE} and 
\eqref{eq:semidiscrete} stochastic. 
As seen in Figure \ref{fig:PDEcontrolcomparison} 
the introduction of this small uncertainty, causes large 
perturbations to the system solution under nominal 
control especially at t=8. To compute the optimal control, 
we minimized the functional \eqref{discretefunct} subject to 
$M=1000$ replicas of the dynamical system \eqref{eq:semidiscrete}. 
Specifically, we considered the single-shooting setting
described in section \ref{sec:multi-shooting} 
with $\Delta t=0.0005$, Euler time integrator, and 
CSE gradient algorithm.  
In Figure \ref{fig:PDEcontrolcomparison} we compare the mean 
and the standard deviation of the solution we obtain by using 
nominal and optimal controls. It is that the optimal control 
is more effective is driving the solution ensemble to zero. 
In fact, both the mean and standard deviation of the 
optimally controlled dynamics are closer to zero. 
Note also that the optimal control differs 
substantially from the nominal control. 

\section{Summary}
\label{sec:summary}
In this paper we developed a new scalable 
algorithm for computational optimal control under uncertainty. 
The new algorithm combines multi-shooting discretization 
methods, interior point optimization, 
common sub-expression elimination, and exact 
gradients obtained via automatic differentiation and 
computational graphs. It also allows for point-wise control constraints
and ensemble state-space constraints. The algorithm has 
an extremely low memory footprint, which allows us to 
process a large number of sample trajectories and 
controls random dynamical systems over long time 
horizons. We also developed a new criterion for verification and 
validation of optimal control based on the stochastic 
version of the Pontryagin's minimum principle. 
We demonstrated the accuracy and computational efficiency 
of the new algorithm in applications to stochastic path planning 
problems involving models of unmanned aerial and ground vehicles, 
and in distributed control of a nonlinear advection-reation-diffusion 
PDE. 

\vs\vs
\noindent 
{\bf Acknowledgements} 
This research was developed with funding from the Defense Advanced Research Projects Agency (DARPA) grant FA8650-18-1-7842.

\bibliographystyle{plain}
\bibliography{bibliography_file}

\end{document}